\documentclass[10pt,reqno]{amsart}
\usepackage{geometry}
\usepackage{hyperref}

\usepackage[utf8]{inputenc}
\usepackage[T1]{fontenc}
\usepackage{amssymb}
\usepackage{mathtools}
\usepackage{enumitem}
\usepackage{xcolor}
\hypersetup{colorlinks=true, linkcolor=blue!60!black, citecolor=blue!60!black, urlcolor=blue!60!black}

\theoremstyle{plain}
\newtheorem{theorem}{Theorem}[section]
\newtheorem{proposition}[theorem]{Proposition}
\newtheorem{lemma}[theorem]{Lemma}
\newtheorem{corollary}[theorem]{Corollary}

\theoremstyle{definition}
\newtheorem{question}[theorem]{Question}
\newtheorem{definition}[theorem]{Definition}
\newtheorem{example}[theorem]{Example}
\newtheorem{remark}[theorem]{Remark}
\newtheorem{openproblem}[theorem]{Open Problem}

\newcommand{\dimH}{\dim_{H}}

\newcommand{\Ha}{\mathcal{H}}
\newcommand{\R}{\mathbb{R}}
\newcommand{\N}{\mathbb{N}}

\newcommand{\Ca}{C_a}

\DeclareMathOperator*{\argmin}{arg\,min}

\begin{document}

\title[Same gaps, different sets]{Same gaps, different sets: a survey on complementary sets and their dimensions}

\author{Nicolás Angelini}
\address[N.~Angelini]{UNSL--IMASL, CONICET, San Luis, Argentina}
\email{nicolas.angelini.2015@gmail.com}

\author{Carlos Cabrelli}
\address[C.~Cabrelli]{Departamento de Matemática, FCEyN, Universidad de Buenos Aires, and
IMAS--CONICET, Buenos Aires, Argentina}
\email{cabrelli@dm.uba.ar}

\author{Ursula Molter}
\address[U.~Molter]{Departamento de Matemática, FCEyN, Universidad de Buenos Aires, and
IMAS--CONICET, Buenos Aires, Argentina}
\email{umolter@dm.uba.ar}

\thanks{The research of the authors was partially supported by Grants PIP 202287/22 (CONICET), and UBACyT 2022-154 (UBA). Nicol\'as Angelini was also partially supported by PROICO 3-0720 “An\'alisis Real y Funcional, Ec.˜Diferenciales.\\ We thank the Fundaci\'on Williams for partially supporting our visit to IMASL} % TODO: financial support (grants), if any

\subjclass[2020]{Primary 28A80; Secondary 28A78, 37C45}

\keywords{Complementary sets, cut-out sets, rearrangements, Cantor sets,
Hausdorff dimension, box dimension, intermediate dimensions, Assouad dimension,
dimension functions, packing measure, thickness}

\begin{abstract}
Let $a=\{a_n\}$ be a non-increasing sequence of positive real numbers adding up
to $1$. Cutting out of the unit interval disjoint open gaps of lengths
$a_1,a_2,\dots$, one for each term, leaves a compact set of Lebesgue measure zero, and
since the gaps can be placed in infinitely many ways, there are infinitely many
such sets. We call them the \emph{complementary sets} of $a$
and write $\mathfrak{C}_a$ for the family of all of them. Every set in
$\mathfrak{C}_a$ is compact and totally disconnected, and all of them
share the same gap lengths; what distinguishes one from another is only the
arrangement of the gaps. Among them, the perfect sets are exactly the Cantor
sets with those gaps, and they form a distinguished subclass.

The central question of this survey is what such a rearrangement can change
and what it cannot. The answer depends sharply on the notion of size: the box
dimensions are constant on $\mathfrak{C}_a$, while the Hausdorff dimension
attains every value in a closed interval whose endpoints are realised by two
distinguished members, the Cantor set $C_a$ built by placing the gaps in a
Cantor-like pattern, and the countable decreasing rearrangement $D_a$. We
survey how this dichotomy persists, and how the roles of $C_a$ and $D_a$ are
exchanged, across the intermediate dimensions of Falconer, Fraser and Kempton,
the Assouad-type and $\Phi$-dimensions, and the Hausdorff and packing
$h$-measures associated with a general dimension function; we also describe the
almost sure behaviour of a randomly rearranged complementary set, where a sharp
threshold in the depth of the scales examined separates two incompatible
regimes. A final part restricts attention to the Cantor subclass and turns to
the dynamical description of the sets $C_a$: the decay conditions under which
$C_a$ is the attractor of an iterated function system, the identification of
the family $C_p=C_{\{n^{-p}\}}$ as an optimally smooth hyperbolic attractor
through Sullivan's scale function, and the consequences for the arithmetic sums
$C_p+C_q$ via thickness and transversality methods.
\end{abstract}

\maketitle

%\setcounter{tocdepth}{1}
%\tableofcontents

% =================================================================
\section{Introduction}
\label{sec:intro}

Start with a non-increasing sequence $a=\{a_n\}$ of positive real numbers
adding up to $1$, and remove from the unit interval $[0,1]$ disjoint open gaps
of lengths $a_1,a_2,\dots$, one for each term of the sequence. Whatever is left
is a compact set of Lebesgue measure zero that contains the endpoints $0$ and
$1$. The sequence, however, does not tell us which set we get: the same lengths
can be laid out along $[0,1]$ in infinitely many ways, and each arrangement
produces a different set. All these sets together form the family
$\mathfrak{C}_a$ of \emph{complementary sets}, or \emph{cut-out sets}, of $a$.
Nothing is lost by this point of view: an infinite compact subset of $[0,1]$ of
measure zero containing the endpoints is determined by the open intervals of its
complement, so it is a complementary set of exactly one such sequence, the
sequence of lengths of those intervals. This family is the object of the
survey. Every member of $\mathfrak{C}_a$ is compact, totally disconnected and
null, but it need not be perfect. The perfect members are exactly the Cantor
sets with gap lengths $a$; they form a distinguished subclass of
$\mathfrak{C}_a$, and we shall restrict to it only when the mathematics requires
it. Two members of $\mathfrak{C}_a$ are canonical. Placing the gaps in a Cantor-like pattern gives a perfect set
$C_a$, the \emph{Cantor set associated with $a$}; placing them in decreasing
order from one endpoint gives a countable set $D_a$ with a single accumulation
point, which is not a Cantor set at all. Both have exactly the same gap lengths,
and yet, as soon as one measures them, they could hardly be more different.

The question that organises this survey is therefore the following: \emph{given
the gap lengths, what is determined for every complementary set in
$\mathfrak{C}_a$, and what is left free to the rearrangement?} The answer
depends entirely on which notion of size one uses, and the resulting picture is
surprisingly varied.

Some quantities are rigid. Falconer \cite{falconer1997techniques} observed that all sets in $\mathfrak{C}_a$
have comparable covering numbers at every scale, so the upper and lower box
dimensions are constant on the whole family: no rearrangement, perfect or not,
can change them. Others vary as much as they possibly can. Besicovitch and
Taylor \cite{BesicovitchTaylor1954} proved that the Hausdorff dimension of the
sets in $\mathfrak{C}_a$ sweeps the full interval $[0,\dim_H C_a]$, with $C_a$
maximal and $D_a$ minimal, and that every intermediate value is attained. This
contrast---one dimension constant, another taking every admissible value---is
the starting point of everything that follows.

Between these two extremes lie the \emph{intermediate dimensions} of Falconer,
Fraser and Kempton, obtained by restricting the diameters of a cover to a
geometric window governed by a parameter $\theta$. Angelini and Molter
\cite{AngeliniMolter2025} computed them for $C_a$ and $D_a$ and showed that these
two sets again bracket the whole family $\mathfrak{C}_a$, with the full interval
attained for every fixed $\theta$: a $\theta$-parametrised version of the
Besicovitch--Taylor theorem, degenerating to it at $\theta=0$ and to the rigid
box dimension at $\theta=1$.

Assouad-type dimensions reverse the picture. They record the extreme local
behaviour of a set rather than its global size, and for the Assouad dimension
and its upper refinements it is $C_a$ that minimises and $D_a$ that maximises
over $\mathfrak{C}_a$, while for the lower Assouad dimension $D_a$ is minimal
and $C_a$ maximal. García, Hare and
Mendivil established this for $\dim_A$ and $\dim_L$
\cite{GarciaHareMendivil2018} and then, through the family of $\Phi$-dimensions
\cite{GarciaHareMendivil2019b}, for a whole one-parameter range interpolating
towards the box dimensions. The same authors \cite{GarciaHareMendivil2019}
identified a sharp threshold in the depth of the scales examined, at
$\log|\log x|/|\log x|$, separating two incompatible almost sure behaviours of a
randomly chosen complementary set in $\mathfrak{C}_a$: below it the Assouad
dimensions are almost surely extremal regardless of $C_a$, above it they almost
surely agree with those of $C_a$.

A further refinement replaces the exponent $s$ by a general gauge $h$. Here the
distinguished Cantor set $C_a$ comes to the fore. Cabrelli, Mendivil, Molter
and Shonkwiler \cite{CabrelliMendivilMolterShonkwiler} showed that every
sequence admits a gauge $h_a$ making $C_a$ an $h_a$-set, even when no power
function does, and García, Molter and Scotto \cite{GarciaMolterScotto2007}
characterised which gauges measure $C_a$ and exhibited Cantor sets with the same
dimension that only a gauge can tell apart. Returning to the whole family, the
reversal appears once more at this finer level: Hare, Mendivil and Zuberman
\cite{HareMendivilZuberman} proved that, up to multiplicative constants, $C_a$
has the maximal Hausdorff $h$-measure and, simultaneously, the minimal packing
$h$-premeasure among all complementary sets in $\mathfrak{C}_a$, and that the
packing measures attained across the family form either $\{0\}$ or all of
$[0,\infty]$, with nothing in between; moreover, prescribed values of these
measures can be realised by perfect sets, that is, within the Cantor subclass.
Not every compact
set, however, is measured by some gauge, and \cite{CabrelliDarjiMolter2011}
shows that this failure is not a rarity.

The recurring theme is that $C_a$ and $D_a$ are extremal in $\mathfrak{C}_a$,
exactly or up to constants, for every notion of size considered here, but which
of the two is the maximum depends on whether the notion is global or local,
covering-based or packing-based. In most cases where a full description of the
attained values is available it takes the same shape: a closed interval with
endpoints at those two sets. The packing measures are the notable exception,
with the all-or-nothing dichotomy described above.

So far the sequence $a$ has been treated as a purely combinatorial datum, and
the results have concerned the whole family $\mathfrak{C}_a$. The sequence can
also be read dynamically, and for this the Cantor subclass, and in fact the
single set $C_a$, is the natural setting. When the gap lengths decay with enough
regularity, the Cantor set $C_a$ turns out to be the attractor of an iterated
function system, and the smoothness of that system is governed precisely by the
rate of decay:
García \cite{GarciaThesis2008} showed that the sets $C_p$ associated with
$a_n=n^{-p}$ are $C^{1+1/p}$-hyperbolic, with $1/p$ the optimal Hölder exponent,
and characterised which central Cantor sets admit realising systems of a given
smoothness. This is not a separate story. The decay condition that makes $C_a$
dynamically realisable is a strengthening of the condition
$0<\mathcal{H}^{s}(C_a)\le\mathcal{P}^{s}_0(C_a)<\infty$ studied in
Section~\ref{sec:complementary}, so the dynamically defined sets are among the
best-behaved members of the theory developed there. In
the other direction, the dynamical description gives information the sequence
alone does not: through Sullivan's scale function, a complete smooth invariant,
$C_p$ is $C^{1+1/p}$-conjugate to a self-similar set, and this conjugacy
transports the whole theory of Bernoulli convolutions to a genuinely non-linear
family.

That transport is what makes the last topic of the survey accessible. Whether
the arithmetic sum $C+C'$ of two Cantor sets on the line contains an interval,
or has positive measure, is a question with a long history in dynamics, going
back to Newhouse's work on homoclinic bifurcations and to the Palis conjecture.
For the family $\{C_p\}$ both the thickness criterion and the transversality
method of Peres and Solomyak can be adapted, giving
$\dim_H(C_p+C_q)=\min\{1/p+1/q,1\}$ for almost every $p$ with $q$ fixed: an
exact non-linear counterpart of the Peres--Shmerkin formula for self-similar
sets.

The survey is organised as follows. Section~\ref{sec:dimensions} collects the
dimensions and measures used throughout: the Hausdorff, box, intermediate and
Assouad-type dimensions, the $\Phi$-dimensions, and dimension functions with
the associated Hausdorff and packing $h$-measures.
Section~\ref{sec:complementary} introduces the family $\mathfrak{C}_a$ of
complementary sets, its Cantor subclass and the distinguished members $C_a$ and
$D_a$, and treats the question of what a rearrangement can and cannot change.
Section~\ref{sec:dynamical} restricts to the Cantor sets $C_a$ and turns to
their dynamical description, and Section~\ref{sec:sums} to sums of Cantor sets
and thickness. The two halves are largely independent: a reader interested only
in rearrangements may stop at the end of Section~\ref{sec:complementary}, and
Sections~\ref{sec:dynamical} and \ref{sec:sums} use from it only the
$p$-Cantor sets of Example~\ref{ex:p-cantor}, Theorem~\ref{thm:hset} and the
Hausdorff and packing $h$-measures of Section~\ref{sec:packing-classification}.
Proofs are
kept to a minimum, and open problems are stated where the picture is incomplete.

% =================================================================
\section{Preliminaries on dimensions}\label{sec:dimensions}

Throughout, $|A|$ denotes the diameter of a non-empty set $A\subset\mathbb{R}^n$,
and we write $a_n\asymp b_n$ to mean that there exist constants $c_1,c_2>0$,
independent of $n$, such that $c_1b_n\le a_n\le c_2b_n$ for all $n$. Given
$\delta>0$, a family $\{U_i\}_{i\in I}$ is a \emph{$\delta$-cover} of $E$ if
$E\subset\bigcup_i U_i$ and $|U_i|\le\delta$ for all $i$. All the sets we
consider from Section~\ref{sec:complementary} onwards are subsets of $[0,1]$,
but the definitions below are stated in $\mathbb{R}^n$ for the sake of
generality.

For $E\subset\mathbb{R}^n$, $s>0$ and $\delta>0$, set
\[
\mathcal{H}^s_\delta(E)=\inf\Big\{\sum_{i=1}^{\infty}|U_i|^s :
\{U_i\}\text{ is a }\delta\text{-cover of }E\Big\}.
\]
The quantity $\mathcal{H}^s_\delta(E)$ is non-decreasing as $\delta$ decreases,
so the limit
\[
\mathcal{H}^s(E)=\lim_{\delta\to0}\mathcal{H}^s_\delta(E)
\]
exists; it is the \emph{$s$-dimensional Hausdorff measure} of $E$. If $t>s$ and
$\{U_i\}$ is a $\delta$-cover of $E$, then
$\sum_i|U_i|^t\le\delta^{t-s}\sum_i|U_i|^s$, whence there is a critical exponent
at which $\mathcal{H}^s(E)$ jumps from $\infty$ to $0$. The \emph{Hausdorff
dimension} of $E$ is
\[
\dim_H E=\inf\{s\ge0:\mathcal{H}^s(E)=0\}
=\sup\{s\ge0:\mathcal{H}^s(E)=\infty\}.
\]
Equivalently, and in the form that will be relevant below,
\begin{equation}\label{eq:hausdorff-cover}
\dim_H E=\inf\Big\{s\ge0:\ \forall\,\varepsilon>0\ \exists\text{ a cover }
\{U_i\}\text{ of }E\text{ with }\sum_i|U_i|^s<\varepsilon\Big\}.
\end{equation}
We refer to \cite{Falconer2014} for the standard
properties of $\mathcal{H}^s$ and $\dim_H$; we shall use repeatedly that
countable sets have Hausdorff dimension $0$.

The Hausdorff measure $\mathcal{H}^s$ assigns a finite positive value to $E$ for
at most one exponent $s$, and for many sets for none at all. Replacing the power
$x^{s}$ by a general gauge $h$ produces a finer scale of measures, and pairing
it with the dual notion obtained from packings rather than covers gives the
framework in which the sets of Section~\ref{sec:complementary} are best
described: the pair refines the classification of Cantor sets by dimension
functions (Section~\ref{sec:packing-classification}) and reveals an extremal
reversal within $\mathfrak{C}_a$ (Section~\ref{sec:rearrangements}) that is
invisible at the level of Hausdorff measure alone.

A \emph{dimension function} is a map $h:(0,\infty)\to\mathbb{R}^{+}$ that is
increasing, continuous, doubling (that is, $h(2x)\le c\,h(x)$ for some constant
$c$) and satisfies $\lim_{x\to0}h(x)=0$. We write $\mathcal{D}$ for the class of
all dimension functions, and $g_s(x)=x^{s}$ for the power functions. 

\begin{definition}\label{def:equiv-dimfun}
Given $g,h\in\mathcal{D}$ we write
\[
g\prec h \iff \lim_{t\to0}\frac{h(t)}{g(t)}=0,
\qquad
g\sim h \iff c_1h(t)\le g(t)\le c_2h(t)\ \text{ for small } t,
\]
and $g\preccurlyeq h$ if $g\prec h$ or $g\sim h$.
\end{definition}

The direction of $\prec$ is worth emphasising: $g\prec h$ means that $h$ is of
\emph{smaller} order than $g$ near the origin, so that $g_s\prec g_t$ exactly
when $s<t$, and the ordering of gauges matches the ordering of the corresponding
dimensions. The same relation controls Hausdorff and packing (pre)measures
simultaneously, which is what allows estimates to be transferred from one to the
other throughout.

Given $E\subset\mathbb{R}$ and $\delta>0$, a \emph{$\delta$-covering} of $E$ is a
countable family of sets of diameter at most $\delta$ whose union contains $E$,
and a \emph{$\delta$-packing} of $E$ is a disjoint family of open intervals of
diameter at most $\delta$ centred at points of $E$. For $h\in\mathcal{D}$ set
\[
\mathcal{H}^h_\delta(E)=\inf\Big\{\sum_{i\ge1}h(|E_i|) :
\{E_i\}\text{ is a }\delta\text{-covering of }E\Big\},
\qquad
\mathcal{H}^h(E)=\lim_{\delta\to0}\mathcal{H}^h_\delta(E),
\]
and
\[
\mathcal{P}^h_\delta(E)=\sup\Big\{\sum_{i\ge1}h(|E_i|) :
\{E_i\}\text{ is a }\delta\text{-packing of }E\Big\},
\qquad
\mathcal{P}^h_0(E)=\lim_{\delta\to0}\mathcal{P}^h_\delta(E).
\]
The set function $\mathcal{H}^h$ is the \emph{Hausdorff $h$-measure}. The
quantity $\mathcal{P}^h_0$, the \emph{$h$-packing premeasure}, is not a measure,
and one obtains the \emph{$h$-packing measure} by the usual regularisation
\[
\mathcal{P}^h(E)=\inf\Big\{\sum_i\mathcal{P}^h_0(E_i) :
E=\bigcup_{i\ge1}E_i\Big\}.
\]
Always $\mathcal{P}^h(E)\le\mathcal{P}^h_0(E)$, and since $h$ is doubling,
$\mathcal{H}^h(E)\le\mathcal{P}^h(E)$. Taking $h=g_s$ recovers the usual
$s$-dimensional Hausdorff and packing (pre)measures, and the index
$\inf\{s:\mathcal{P}^{g_s}_0(E)=0\}$---the \emph{pre-packing dimension} of $E$---
coincides with $\overline{\dim}_B E$. The gap between $\mathcal{P}^h_0$ and
$\mathcal{P}^h$ is not a technicality here: a countable set has zero packing
measure for every $h$, but may well have infinite premeasure.

Finally, $E$ is called an \emph{$h$-set} if
$0<\mathcal{H}^h(E)<\infty$, and $h$ is \emph{subadditive} if
$h(x+y)\le C(h(x)+h(y))$ for all $x,y$ and some constant $C$. 

Not every set is an $h$-set; the simplest example is countable: if $E=\{1/n:n\in\mathbb{N}\}\cup\{0\}$, then
$\mathcal{H}^h$ is non-atomic for every $h\in\mathcal{D}$, since $h(x)\to0$ as
$x\to0$, so countable additivity forces $\mathcal{H}^h(E)=0$ for every gauge.
This example is unsatisfying, however, because $E$ is countable and hence of
Hausdorff dimension $0$ for trivial reasons. A genuinely fractal example is due
to Best \cite{Best1939}.

The same dichotomy appears for sets arising in Diophantine approximation. Olsen and Renfro
\cite{OlsenRenfro2006} located the exact cut point between the gauges
for which the set $\mathbb{L}$ of Liouville numbers has infinite (indeed
non-$\sigma$-finite) Hausdorff $h$-measure and those for which it has measure
zero, so $\mathbb{L}$ is not dimensional; more recently, Moreira, Erazo and
Angelini \cite{MoreiraErazoAngelini} obtained the analogous description for the
irrationals whose Lagrange value equals $3$, and showed that $\mathbb{L}$ is
strictly the larger of the two on this finer scale: there are gauges $h$ with
$\mathcal{H}^{h}(\mathbb{L})=\infty$ and $\mathcal{H}^{h}(k^{-1}(3))=0$.

Cabrelli, Darji and Molter \cite{CabrelliDarjiMolter2011} placed this question
in a broader setting by allowing measures other than Hausdorff measures. A
compact set $M$ in a Polish group is \emph{$\mathcal{H}$-visible} if it is an
$h$-set for some $h\in\mathcal{D}$, and \emph{visible} if $0<\mu(M)<\infty$ for
some translation-invariant Borel measure $\mu$, not necessarily of Hausdorff
type; $M$ is \emph{strongly invisible} if every translation-invariant Borel
measure either vanishes on $M$ or fails to be $\sigma$-finite there. Not being an $h$-set for any $h$ is thus weaker than being strongly invisible, and Best's example
is not known to be strongly invisible.
Davies had already exhibited a strongly invisible compact subset of $\mathbb{R}$ (see \cite{Davies1971}); Cabrelli, Darji and Molter extend the construction to every uncountable abelian Polish group.

For a bounded set $E\subset\mathbb{R}^n$, let $N_\delta(E)$ denote the smallest
number of closed balls of radius $\delta$ needed to cover $E$. The \emph{lower}
and \emph{upper box dimensions} of $E$ are
\[
\underline{\dim}_B E=\liminf_{\delta\to0}\frac{\log N_\delta(E)}{-\log\delta},
\qquad
\overline{\dim}_B E=\limsup_{\delta\to0}\frac{\log N_\delta(E)}{-\log\delta},
\]
and when both coincide we write $\dim_B E$ for the common value. Comparing
covers by balls, by grid cubes of side $\delta$ and by arbitrary sets of
diameter at most $\delta$ changes $N_\delta(E)$ only by a multiplicative
constant and a bounded rescaling of $\delta$, so $N_\delta(E)$ may equally be
taken to be the number of grid cubes of side $\delta$ meeting $E$, the least
number of sets of diameter at most $\delta$ covering $E$, or the largest number
of disjoint balls of radius $\delta$ centred in $E$.

For bounded $E$, the box dimensions admit a reformulation parallel to
\eqref{eq:hausdorff-cover}, in which the exponent is again characterised through
covering sums rather than counting functions:
\begin{equation}\label{eq:lower-box-cover}
\underline{\dim}_B E=\inf\left\{
\begin{array}{ll}
s\ge0: & \forall\,\varepsilon>0\ \forall\,\delta_0>0\ \exists\,\delta\in(0,\delta_0)\\
& \text{and a cover }\{U_i\}\text{ of }E\text{ with }|U_i|=\delta\ \forall i\\
& \text{such that }\sum_i|U_i|^s<\varepsilon
\end{array}\right\}
\end{equation}
and
\begin{equation}\label{eq:upper-box-cover}
\overline{\dim}_B E=\inf\left\{
\begin{array}{ll}
s\ge0: & \forall\,\varepsilon>0\ \exists\,\delta_0>0\text{ such that }
\forall\,\delta\in(0,\delta_0)\\
& \text{there is a cover }\{U_i\}\text{ of }E\text{ with }|U_i|=\delta\ \forall i\\
& \text{such that }\sum_i|U_i|^s<\varepsilon
\end{array}\right\}.
\end{equation}

Two further dimensions will be needed, both of which describe the extreme local
behaviour of a set rather than its global size. Given a non-empty
$E\subset\mathbb{R}^n$, recall that $N_r(E)$ denotes the smallest number of
closed balls of radius $r$ needed to cover $E$. The \emph{Assouad dimension} of
$E$, is defined as follows
\begin{align*}
\dim_A E=\inf\Big\{\alpha :\ &\exists\,c,\rho>0\text{ such that }\\
&\sup_{x\in E}N_r\big(B(x,R)\cap E\big)\le c\Big(\frac{R}{r}\Big)^{\alpha}
\text{ for all }0<r<R<\rho\Big\},
\end{align*}
and its dual, the \emph{lower Assouad dimension} is
\begin{align*}
\dim_L E=\sup\Big\{\alpha :\ &\exists\,c,\rho>0\text{ such that }\\
&\inf_{x\in E}N_r\big(B(x,R)\cap E\big)\ge c\Big(\frac{R}{r}\Big)^{\alpha}
\text{ for all }0<r<R<\rho\Big\}.
\end{align*}
Whereas the Hausdorff dimension records global geometric information, these two
quantities record the coarser but local behaviour of $E$, since only one ball
$B(x,R)$ is examined at a time. They are correspondingly sensitive to local
irregularities: if $E$ has an isolated point, then $\dim_L E=0$. For a compact
set $E$ they bracket all the dimensions considered above,
\begin{equation}\label{eq:chain}
\dim_L E\le\dim_H E\le\underline{\dim}_B E
\le\overline{\dim}_B E\le\dim_A E.
\end{equation}
For further information on Assouad-like dimensions we refer the reader to
Fraser's book \cite{fraser2020assouad}.

Hausdorff and box dimension formulations \eqref{eq:hausdorff-cover} and \eqref{eq:lower-box-cover}, differ in
a single structural respect: the Hausdorff dimension places no constraint on the
diameters of the covering sets, so that countable covers by sets of widely
differing sizes are admissible, whereas the box dimensions force all the sets of
the cover to have the same diameter $\delta$, and hence the cover to be finite.
Since the class of admissible covers is larger in \eqref{eq:hausdorff-cover},
the infimum is smaller, and one obtains at once
\[
\dim_H E\le\underline{\dim}_B E\le\overline{\dim}_B E ,
\]
inequalities which may be strict.

The two definitions therefore appear as the extreme cases of a single scheme.
This is the point of view adopted by Falconer, Fraser and Kempton
\cite{FalconerFraserKempton}, who interpolated between them by allowing the
diameters of the covering sets to range over a geometric window governed by a
parameter $\theta$.

\begin{definition}\label{def:intermediate}
Let $E\subset\mathbb{R}^n$ be bounded and let $0\le\theta\le1$. The \emph{lower
$\theta$-intermediate dimension} of $E$ is
\begin{align*}
\underline{\dim}_\theta E=\inf\Big\{ s\ge0 :\ &\forall\,\varepsilon>0\
\forall\,\delta_0>0\ \exists\,0<\delta\le\delta_0\text{ and a cover }
\{U_i\}\text{ of }E\\
&\text{with }\delta^{1/\theta}\le|U_i|\le\delta\text{ and }
\sum_i|U_i|^s\le\varepsilon\Big\},
\end{align*}
and the \emph{upper $\theta$-intermediate dimension} of $E$ is
\begin{align*}
\overline{\dim}_\theta E=\inf\Big\{ s\ge0 :\ &\forall\,\varepsilon>0\
\exists\,\delta_0>0\text{ such that }\forall\,0<\delta\le\delta_0\\
&\text{there is a cover }\{U_i\}\text{ of }E\\
&\text{with }\delta^{1/\theta}\le|U_i|\le\delta\text{ and }
\sum_i|U_i|^s\le\varepsilon\Big\}.
\end{align*}
When the two values agree we write $\dim_\theta E$ for the common value.
\end{definition}

With the convention $\delta^{1/0}=0$, Definition~\ref{def:intermediate} recovers
the two extremes:
\[
\underline{\dim}_0E=\overline{\dim}_0E=\dim_H E,\qquad
\underline{\dim}_1E=\underline{\dim}_BE,\qquad
\overline{\dim}_1E=\overline{\dim}_BE .
\]
For $\theta\in[0,1]$ the admissible covers are more numerous than those allowed
in \eqref{eq:lower-box-cover} and fewer than those allowed in
\eqref{eq:hausdorff-cover}, so that
\[
\dim_H E\le\underline{\dim}_\theta E\le\overline{\dim}_\theta E
\le\overline{\dim}_B E
\qquad\text{and}\qquad
\underline{\dim}_\theta E\le\underline{\dim}_B E ,
\]
and, by the same reasoning, $\theta\mapsto\underline{\dim}_\theta E$ and
$\theta\mapsto\overline{\dim}_\theta E$ are non-decreasing, though not
necessarily strictly increasing.

We record the properties that will be used later. Both dimensions are monotone
under inclusion and invariant under bi-Lipschitz maps, and satisfy
$\underline{\dim}_\theta E=\underline{\dim}_\theta\overline{E}$ and
$\overline{\dim}_\theta E=\overline{\dim}_\theta\overline{E}$ for
$\theta\in(0,1]$. The upper intermediate dimension is finitely stable,
\[
\overline{\dim}_\theta(E\cup E')=
\max\{\overline{\dim}_\theta E,\overline{\dim}_\theta E'\},
\]
whereas $\underline{\dim}_\theta$ is not finitely stable, and neither dimension
is countably stable. It was shown in \cite{FalconerFraserKempton} that
$\theta\mapsto\underline{\dim}_\theta E$ and
$\theta\mapsto\overline{\dim}_\theta E$ are continuous on $(0,1]$; continuity at
$\theta=0$, however, may fail, so that the intermediate dimensions need not
interpolate all the way down to $\dim_H E$.

These dimensions have proved to be more than an interpolation device. Being a
function rather than a number, $\theta\mapsto\dim_\theta E$ is a finer
bi-Lipschitz invariant than its endpoints: Banaji and Kolossváry
\cite{BanajiKolossvary2024} computed it exactly for Bedford--McMullen carpets
and used it to obtain necessary conditions for Lipschitz equivalence. Its
behaviour as $\theta$ tends to $0$ also characterises when the box dimension of
a typical orthogonal projection attains its maximal value
\cite{AngeliniMolterProjections}. Knowing the whole spectrum therefore says more than
knowing its two endpoints.

\begin{example}\label{ex:Fp}
Let $p>0$ and
$F_p=\{0\}\cup\{n^{-p}:n\in\mathbb{N}\}$. Since $F_p$ is countable,
$\dim_H F_p=0$, while $\dim_B F_p=1/(p+1)$. It was proved in
\cite{FalconerFraserKempton} that
\[
\underline{\dim}_\theta F_p=\overline{\dim}_\theta F_p=\frac{\theta}{p+\theta},
\qquad \theta\in(0,1] .
\]
Here the interpolation is complete: the dimension function is continuous on all
of $[0,1]$ and tends to $\dim_H F_p=0$ as $\theta\to0^+$. In the language of
Section~\ref{sec:complementary}, $F_p$ is the decreasing rearrangement $D_a$
of the sequence $a_n=n^{-p}-(n+1)^{-p}$, and the formula above is the case
$\overline{\dim}_B C_a=1/(p+1)$ of Theorem~\ref{thm:extremos}\,(ii).
\end{example}

Beyond these general properties, one may ask which functions of $\theta$ occur
as intermediate dimensions. The answer, obtained by Banaji and Rutar
\cite{BanajiRutar2022}, is a complete characterisation, and the constraint turns
out to be governed by the two Assouad dimensions of the set.

The starting point is a quantitative bound: the intermediate dimensions cannot
grow arbitrarily fast in $\theta$.

\begin{proposition}[{\cite[Cor.~2.8]{BanajiRutar2022}}]\label{prop:Banaji-bound}
If $\dim_L E<\dim_A E$, then for every $\theta\in(0,1]$,
\[
\overline{\dim}_\theta E\ \ge\
\frac{\theta\dim_A E\,(\overline{\dim}_B E-\dim_L E)
+\dim_L E\,(\dim_A E-\overline{\dim}_B E)}
{\theta\,(\overline{\dim}_B E-\dim_L E)+(\dim_A E-\overline{\dim}_B E)} ,
\]
and the same holds with $\overline{\dim}_\theta$ and $\overline{\dim}_B$
replaced by $\underline{\dim}_\theta$ and $\underline{\dim}_B$.
\end{proposition}

Writing $\lambda=\dim_L E$, $\alpha=\dim_A E$ and $\beta=\overline{\dim}_B E$,
the right-hand side is exactly the solution of
\[
h'(\theta)=\frac{(h(\theta)-\lambda)(\alpha-h(\theta))}{(\alpha-\lambda)\,\theta},
\qquad h(1)=\beta.
\]
 A dimension profile therefore cannot
decrease faster, as $\theta$ decreases, than the extremal solution of this
equation. This makes it natural to single out the functions satisfying the
corresponding differential inequality.

\begin{definition}\label{def:BRclass}
For $0\le\lambda\le\alpha\le n$, let $\mathfrak{B}(\lambda,\alpha)$ be the class
of functions $h:[0,1]\to[\lambda,\alpha]$ that are non-decreasing, continuous on
$(0,1]$, and satisfy
\[
\limsup_{\epsilon\to0^{+}}\frac{h(\theta+\epsilon)-h(\theta)}{\epsilon}
\le\frac{(h(\theta)-\lambda)(\alpha-h(\theta))}{(\alpha-\lambda)\,\theta}
\qquad\text{for all }\theta\in(0,1).
\]
If $\lambda=\alpha$ we set $\mathfrak{B}(\lambda,\lambda)=\{h\equiv\lambda\}$.
\end{definition}

\begin{theorem}[{\cite{BanajiRutar2022}}]\label{thm:banajirutar}
Let $E\subset\mathbb{R}^n$ be bounded with $\dim_L E=\lambda$ and
$\dim_A E=\alpha$. Then the functions $h_1(\theta)=\underline{\dim}_\theta E$
and $h_2(\theta)=\overline{\dim}_\theta E$ satisfy
\[
h_1,h_2\in\mathfrak{B}(\lambda,\alpha),\qquad h_1\le h_2,\qquad h_1(0)=h_2(0).
\]
Conversely, given $h_1,h_2\in\mathfrak{B}(\lambda,\alpha)$ with $h_1\le h_2$ and
$h_1(0)=h_2(0)$, there is a compact perfect set $E\subset\mathbb{R}^n$ with
\[
\dim_L E=\lambda,\quad \underline{\dim}_\theta E=h_1(\theta),\quad
\overline{\dim}_\theta E=h_2(\theta),\quad \dim_A E=\alpha .
\]
\end{theorem}

Theorem~\ref{thm:banajirutar} settles the question for arbitrary compact sets.
The corresponding question inside a fixed family $\mathfrak{C}_a$ is open, and we
return to it at the end of Section~\ref{sec:complementary}.

Finally, we mention two extensions of this circle of ideas. In
\cite{Banaji2023}, Banaji introduced the $\Phi$-intermediate dimensions, in
which the admissible diameters are constrained to $[\Phi(\delta),\delta]$,
recovering the $\theta$-intermediate dimensions for $\Phi(x)=x^{1/\theta}$ and
extending the theory to general metric spaces and in
\cite{AngeliniMolterTejada2026}, intermediate dimensions were defined for finite
Borel measures, interpolating between the upper Hausdorff and the Minkowski
dimensions of a measure.

Just as the intermediate dimensions interpolate between $\dim_H$ and $\dim_B$,
Lü and Xi \cite{lu2016quasi} introduced a spectrum interpolating between the
upper box dimension and the Assouad dimension. For non-empty
$E\subset\mathbb{R}^n$ and $\theta\in(0,1]$,
\begin{align*}
\dim_A^\theta E=\inf\Big\{\alpha :\ &\exists\,C>0\text{ such that }\\
&\sup_{x\in E}N_r\big(B(x,R)\cap E\big)\le C\Big(\frac{R}{r}\Big)^{\alpha}\\
&\text{for all }0<r\le R^{1/\theta}<R<1\Big\},
\end{align*}
so that $\dim_A^1E=\dim_A E$. The \emph{quasi-Assouad dimension} is then
defined as $\dim_{qA}E=\lim_{\theta\to1}\dim_A^\theta E$. In general
\[
\overline{\dim}_B E\le\dim_A^\theta E\le\dim_{qA}E\le\dim_A E ,
\]
and all the inequalities may be strict.

The Assouad dimensions admit refinements that interpolate towards the box
dimensions, and these are of a different nature from the intermediate
dimensions. Whereas $\underline{\dim}_\theta$ and $\overline{\dim}_\theta$
restrict the global range of diameters admitted in a cover, the refinements
below localise around each point and each scale, in the spirit of $\dim_A$
itself. The most general such family is the one introduced by García, Hare and
Mendivil \cite{GarciaHareMendivil2019b, GarciaHareMendivil2019}.

\begin{definition}\label{def:phi-dimensions}
An \emph{admissible function} is a map $\Phi:(0,1)\to\mathbb{R}^{+}$ such that
$R^{1+\Phi(R)}$ decreases as $R\downarrow0$; in particular
$R^{1+\Phi(R)}\le R$. Given a bounded, non-empty $E\subset\mathbb{R}^n$, the
\emph{upper $\Phi$-dimension} of $E$ is
\begin{align*}
\overline{\dim}_\Phi E=\inf\Big\{\alpha :\ &\exists\,c_1,c_2>0\text{ such that }\\
&\sup_{x\in E}N_r\big(B(x,R)\cap E\big)\le c_2\Big(\frac{R}{r}\Big)^{\alpha}\\
&\text{for all }0<r\le R^{1+\Phi(R)}\le R<c_1\Big\},
\end{align*}
and the \emph{lower $\Phi$-dimension} is
\begin{align*}
\underline{\dim}_\Phi E=\sup\Big\{\alpha :\ &\exists\,c_1,c_2>0\text{ such that }\\
&\inf_{x\in E}N_r\big(B(x,R)\cap E\big)\ge c_2\Big(\frac{R}{r}\Big)^{\alpha}\\
&\text{for all }0<r\le R^{1+\Phi(R)}\le R<c_1\Big\}.
\end{align*}
\end{definition}

The choice $\Phi\equiv0$ recovers $\dim_A E$ and $\dim_L E$. The constant
function $\Phi\equiv1/\theta-1$, $\theta\in(0,1)$, gives the Assouad spectrum  of Lü and Xi, for which we write
$\dim_A^\theta E=\overline{\dim}_\Phi E$; letting $\Phi(R)\to0$ as $R\to0$
yields the \emph{quasi-Assouad dimensions},
\[
\dim_{qA}E=\lim_{\theta\to1}\dim_A^\theta E,
\]
and their lower counterpart $\dim_{qL}E$. For compact $E$ one has
\begin{align*}
\dim_L E &\le\dim_{qL}E\le\dim_H E\le\underline{\dim}_B E
\le\overline{\dim}_B E\le\dim_{qA}E\le\dim_A E,\\
\dim_L E &\le\underline{\dim}_\Phi E\le\underline{\dim}_B E
\le\overline{\dim}_B E\le\overline{\dim}_\Phi E\le\dim_A E,
\end{align*}
and all the inequalities may be strict.

% =================================================================

% =================================================================
\section{Cut-out sets and their rearrangements}
\label{sec:complementary}
% =================================================================

Let $a=\{a_n\}$ be a positive, non-increasing, summable sequence with sum $1$,
and let $I_a=[0,1]$. We denote by $\mathfrak{C}_a$ the family of all compact sets
$E\subset I_a$ of the form $E=I_a\setminus\bigcup_{j\ge1}U_j$, where $\{U_j\}$
is a pairwise disjoint family of open intervals contained in $I_a$ with
$|U_j|=a_j$ for every $j$. The sets in $\mathfrak{C}_a$ are called the
\emph{complementary sets} of $a$. Every $E\in\mathfrak{C}_a$ is compact, totally
disconnected and of Lebesgue measure zero, and contains the endpoints $0$ and
$1$; conversely, every infinite compact subset of $[0,1]$ of measure zero
containing $0$ and $1$ belongs to $\mathfrak{C}_a$ for a unique such sequence $a$. Note that
$\mathfrak{C}_a$ depends only on the collection of lengths $\{a_j\}$ and not on
the order in which they are listed, so no generality is lost in assuming $a$
non-increasing.

Two members of $\mathfrak{C}_a$ will play a distinguished role. The first is the
countable set
\[
D_a=\{0\}\cup\Big\{\sum_{i\ge k}a_i : k\in\mathbb{N}\Big\},
\]
a decreasing sequence of points accumulating at $0$, whose consecutive gaps are
exactly the terms of $a$ arranged in decreasing order. We call $D_a$ the
\emph{decreasing rearrangement} of $a$.

The second is the \emph{Cantor set associated with $a$}, denoted $C_a$ and
constructed inductively as follows. In the first step, remove from $I_a$ an open
interval of length $a_1$, leaving two closed intervals $I_{1,1}$ and $I_{1,2}$.
Having completed the $k$-th step, we have $2^{k}$ closed intervals
$I_{k,1},\dots,I_{k,2^{k}}$; from each $I_{k,j}$ we remove an open interval of
length $a_{2^{k}+j-1}$, obtaining the two closed intervals $I_{k+1,2j-1}$ and
$I_{k+1,2j}$. We then set
\[
C_a:=\bigcap_{k\ge1}\bigcup_{j=1}^{2^{k}}I_{k,j}.
\]
The construction determines $C_a$ uniquely, since the lengths of the gaps
removed on either side of a given gap are prescribed; for instance, the interval
of length $a_1$ must be positioned so that the remaining closed interval on its
left has length $a_2+a_4+a_5+a_8+\cdots$. The classical middle-third Cantor set
is $C_a$ for the sequence given by $a_j=3^{-n}$ when $2^{n-1}\le j<2^{n}$.

Since every $E\in\mathfrak{C}_a$ is already compact, totally disconnected and of
measure zero, such a set is homeomorphic to the classical Cantor set precisely
when it is perfect. We shall accordingly refer to a perfect $E\in\mathfrak{C}_a$
as a \emph{Cantor rearrangement} of $a$, reserving the notation $C_a$ for the
particular set constructed above. Not every complementary set is of this type:
$D_a$ consists entirely of isolated points except for $0$.

Given a non-increasing, summable sequence $a=\{a_n\}$, define
\[
\displaystyle x_n(a)=\sum_{i=n}^{\infty}a_i
\quad \text{ and }\quad s_n(a)=2^{-n} x_{2^n}(a).\]

When the sequence $a$ is fixed and clear from context, we will simply write $s_n$ instead of $s_n (a)$.

Note that since $a$ is non increasing we have the estimate  
\[s_{n+1} < |I_{n,j}| < s_{n-1}
\quad \text{for all $1 \leq j \leq 2^n$,}\]
where $I_{n,j}$ are the closed intervals of step $n$ in the construction of the Cantor set $C_{a}$.

\subsection{Hausdorff measure and dimension}
The Hausdorff dimension of the sets in $\mathfrak{C}_a$ was determined by
Besicovitch and Taylor \cite{BesicovitchTaylor1954}, who showed that $C_a$ is
extremal in the family and that every smaller value is attained. Throughout we
write $\alpha(a):=\dim_H C_a$.

\begin{theorem}[{\cite{BesicovitchTaylor1954}}]\label{thm:BT}
Let $a=\{a_n\}$ be a positive, non-increasing, summable sequence. Then
\[
\alpha(a)=\liminf_{n\to\infty}\frac{\log n}{-\log(x_n/n)} ,
\]
and $\dim_H E\le\alpha(a)$ for every $E\in\mathfrak{C}_a$. Conversely, for every
$\beta\in[0,\alpha(a)]$ there is a set $E\in\mathfrak{C}_a$ with
$\dim_H E=\beta$. If moreover $\alpha(a)>0$, the measure itself can be
prescribed: for every $\beta\in[0,\alpha(a)]$ and every
$0\le\gamma\le\Lambda(\beta,a)$ there is $E\in\mathfrak{C}_a$ with
$\dim_H E=\beta$ and $\mathcal{H}^{\beta}(E)=\gamma$, where
\[
\Lambda(\beta,a):=\liminf_{n\to\infty}n\Big(\frac{x_n}{n}\Big)^{\beta}.
\]
\end{theorem}

Thus $\{\dim_H E : E\in\mathfrak{C}_a\}=[0,\alpha(a)]$, the minimum being
attained by the decreasing rearrangement $D_a$, which is countable and hence has
Hausdorff dimension $0$, and the maximum by $C_a$. The quantity
$\Lambda(\beta,a)$ reappears in Section~\ref{sec:rearrangements} as the case
$h=g_\beta$ of the estimates of García, Molter and Scotto for $\mathcal{H}^h(C_a)$.

It was shown in \cite{GarciaMolterScotto2007} that the two lower
dimensions agree on $C_a$,
\[
\dim_H C_a=\underline{\dim}_B C_a ,
\]
which is what makes the interpolation degenerate from below in
Theorem~\ref{thm:extremos}\,(i).

\begin{example}[$p$-Cantor sets]\label{ex:p-cantor}
Let $p>1$ and let $\lambda=\{\lambda_k\}$ with $\lambda_k=k^{-p}$, normalised so
that $\sum_k\lambda_k=1$. Then $x_n\asymp n^{1-p}$ and
Theorem~\ref{thm:BT} gives $\dim_H C_\lambda=1/p$. Consequently
$\dim_H E\le1/p$ for every $E\in\mathfrak{C}_\lambda$, and for each
$s\in[0,1/p]$ there is a rearrangement $E\in\mathfrak{C}_\lambda$ with
$\dim_H E=s$ and $\mathcal{H}^{s}(E)>0$.
\end{example}

\subsection{Dimensionality and visibility}
\label{sec:visibility}

Passing from a single exponent $s$ to a general gauge $h$ enlarges considerably
the class of sets that can be assigned a finite positive size, but it does not
exhaust it. Call a compact set $E$ \emph{dimensional} if $E$ is an $h$-set for
some $h\in\mathcal{D}$. 
Not every Cantor set $C_a$ is an $s$-set. What is true, and this is the point of the
following theorem, is that a suitable gauge always exists.

\begin{theorem}[{\cite{CabrelliMendivilMolterShonkwiler}}]\label{thm:hset}
Let $a=\{a_k\}$ be positive, non-increasing, with $\sum_k a_k=1$, and let $C_a$
be the associated Cantor set. Then there is a dimension function $h\in\mathcal{D}$
for which $C_a$ is an $h$-set; in fact $h$ may be chosen concave, and then
\[
\tfrac{1}{32}\le\mathcal{H}^{h}(C_a)\le1 .
\]
\end{theorem}

The gauge function is constructed explicitly from the sequence.

Any $h\in\mathcal{D}$ with $h(x_n/n)\sim 1/n$ is called an \emph{associated
dimension function} of $a$, and we write $h_a$ for any such gauge; all of them
are comparable, so $h_a$ is well defined up to $\sim$.

\subsection{Packing measures and the classification of gauges}
\label{sec:packing-classification}

The two orderings on $\mathcal{D}$ translate directly into statements about
measures, and the translation is the same for Hausdorff measure, packing
premeasure and packing measure.

\begin{lemma}[{\cite{GarciaMolterScotto2007}}]\label{lem:comparison}
Let $\nu^{g}$ denote any of $\mathcal{H}^{g}$, $\mathcal{P}^{g}_0$ or
$\mathcal{P}^{g}$, for $g\in\mathcal{D}$, and let $E\subset\mathbb{R}$.
\begin{enumerate}[label=\textup{(\roman*)}]
\item If $f\prec h$ and $\nu^{f}(E)<\infty$, then $\nu^{h}(E)=0$.
\item If $f\sim h$, then $\nu^{f}(E)<\infty$ if and only if $\nu^{h}(E)<\infty$,
and $\nu^{f}(E)>0$ if and only if $\nu^{h}(E)>0$. In particular $E$ is an
$f$-set precisely when it is an $h$-set, and likewise for the packing analogues.
\end{enumerate}
\end{lemma}

The Besicovitch--Taylor estimates for $\mathcal{H}^{s}(C_a)$ extend to arbitrary
gauges, and admit an exact packing dual.

\begin{theorem}[{\cite[Prop.~4.1, Thm.~4.2]{GarciaMolterScotto2007}}]
\label{thm:gms-packing}
Let $a$ be positive, non-increasing and summable, and put $b_n:=x_n/n$. Then for
every $h\in\mathcal{D}$,
\[
\tfrac14\liminf_{n\to\infty}n\,h(b_n)\le\mathcal{H}^{h}(C_a)
\le 4\liminf_{n\to\infty}n\,h(b_n),
\]
\[
\tfrac18\limsup_{n\to\infty}n\,h(b_n)\le\mathcal{P}^{h}_0(C_a)
\le 8\limsup_{n\to\infty}n\,h(b_n).
\]
\end{theorem}

The two estimates differ only in the replacement of $\liminf$ by $\limsup$, and
of covering by packing. Since an associated dimension function satisfies
$h_a(b_n)\sim1/n$, both quantities $n\,h_a(b_n)$ are bounded above and below, so
Theorem~\ref{thm:gms-packing} gives at once that $\mathcal{H}^{h_a}(C_a)$ and
$\mathcal{P}^{h_a}_0(C_a)$ are finite and positive: the same gauge regularises
the covering and the packing process simultaneously.

Combining the two halves characterises the gauges that measure $C_a$, and in
particular tells when the associated dimension function reduces to a power.

\begin{theorem}[{\cite{GarciaMolterScotto2007}}]\label{thm:gms-characterisation}
Let $a$ be non-increasing and let $h\in\mathcal{D}$ satisfy
$0<\mathcal{H}^{h}(C_a)\le\mathcal{P}^{h}_0(C_a)<\infty$. Then, for every
$g\in\mathcal{D}$,
\[
g\sim h\iff 0<\mathcal{H}^{g}(C_a)\le\mathcal{P}^{g}_0(C_a)<\infty,\]
and
\[
g\prec h\iff\mathcal{H}^{g}(C_a)=+\infty,
\qquad
h\prec g\iff\mathcal{P}^{g}_0(C_a)=0 .
\]
In particular $h\sim g_s$ if and only if
$0<\mathcal{H}^{s}(C_a)\le\mathcal{P}^{s}_0(C_a)<\infty$.
\end{theorem}

Being an $s$-set is therefore not enough for the associated dimension function
to be equivalent to a power: the condition $0<\mathcal{H}^{s}(C_a)<\infty$ must
be supplemented by $\mathcal{P}^{s}_0(C_a)<\infty$. That the extra condition is
not automatic is shown by the following.

\begin{example}[{\cite{GarciaMolterScotto2007}}]\label{ex:gms-counterexample}
For every $s\in(0,1)$ there is a non-increasing sequence $a$, obtained by
perturbing a geometric sequence as $\lambda_k=2^{-k/(s+\varepsilon_k)}$ with
$\{\varepsilon_k\}$ a null sequence supported on alternate dyadic blocks, such
that
\[
\dim_H C_a=\underline{\dim}_B C_a=\overline{\dim}_B C_a=s,
\qquad 0<\mathcal{H}^{s}(C_a)<\infty,
\qquad \mathcal{P}^{s}_0(C_a)=+\infty .
\]
Thus $C_a$ is an $s$-set whose box dimension exists, and yet $h_a\not\sim g_s$.
\end{example}

The gauges therefore classify these sets more finely than the exponents alone.
For $p>1$ and $q\in\mathbb{R}$, the sequence $a_n=(\log n)^{q}n^{-p}$, $n>1$,
produces Cantor sets $C_{p,q}$ with $\dim_H C_{p,q}=\dim_B C_{p,q}=1/p$
irrespective of $q$, while
\[
h_{p,q}(x):=x^{1/p}\,(-\log x)^{-q/p}
\]
satisfies $0<\mathcal{H}^{h_{p,q}}(C_{p,q})\le
\mathcal{P}^{h_{p,q}}_0(C_{p,q})<\infty$, see \cite[Corollary 4.8]{GarciaMolterScotto2007} . Ordering the pairs $(1/p,q)$
lexicographically, the family $\{h_{p,q}\}$ separates sets that the power
functions $\{g_s\}$ cannot tell apart.

\subsection{Hausdorff and packing $h$-measures of rearrangements}
\label{sec:rearrangements}

The results of Section~\ref{sec:packing-classification} concern the single set
$C_a$. We now turn to the whole family, where the natural question is how
$\mathcal{H}^h$ and $\mathcal{P}^h$ vary across $\mathfrak{C}_a$ for a fixed
gauge. Hare, Mendivil and Zuberman \cite{HareMendivilZuberman} extended both the
Besicovitch--Taylor extremality statement and its packing dual to arbitrary
$h\in\mathcal{D}$, and found that the two behave in opposite ways.

The statements below are sharper than mere attainability of values: they locate
the realising sets inside $\mathfrak{C}_a$.

\begin{proposition}[\cite{HareMendivilZuberman}]
\label{prop:hmz-hausdorff-max}
For $h\in\mathcal D$, there is a constant $c$, depending only on the doubling constant of $h$, such that $\Ha^h(E) \le c\,\Ha^h(\Ca)$ for every $E\in\mathfrak{C}_a$. That is, $\Ca$ has, up to a constant, the \emph{maximal} Hausdorff $h$-measure within $\mathfrak{C}_a$, generalising the Besicovitch--Taylor extremality result beyond power functions $g_s$. When $\Ca$ is a middle-$\tau$ Cantor set the constant may be taken to be $c=1$, so $\Ca$ is then \emph{exactly} maximal; in general the sharp constant and the identity of a maximiser (if any exists) remain open.
\end{proposition}

\begin{theorem}[\cite{HareMendivilZuberman}]
\label{thm:hmz-hausdorff-attain}
Let $h\prec h_a$, where $h_a$ is the associated dimension function of $\Ca$, and let $\gamma\ge0$. There is a \emph{perfect} set $E\in\mathfrak{C}_a$ with $\Ha^h(E)=\gamma$. Moreover the set $\{E\in\mathfrak{C}_a : \Ha^h(E)=\gamma\}$ is dense in $\mathfrak{C}_a$ with respect to the Hausdorff metric $\rho$.
\end{theorem}

Theorem~\ref{thm:hmz-hausdorff-attain} sharpens the Besicovitch--Taylor attainability statement---that any $\beta\in[0,\alpha(a)]$ is realised as $\dimH E$ for some $E\in\mathfrak{C}_a$---to the level of $h$-measures rather than merely dimension values, and it guarantees that the realising set can always be chosen perfect, and in fact densely so, within $(\mathfrak{C}_a,\rho)$.

The dual statement for packing measure exhibits the \emph{opposite} extremal pattern, dual to Proposition~\ref{prop:hmz-hausdorff-max}.

\begin{proposition}[\cite{HareMendivilZuberman}]
\label{prop:hmz-packing-min}
There is a constant $c$ such that $\mathcal P^h_0(\Ca) \le c\, \mathcal P^h_0(E)$ for every $E\in\mathfrak{C}_a$ and every $h\in\mathcal D$. That is, despite having the maximal Hausdorff $h$-measure by Proposition~\ref{prop:hmz-hausdorff-max}, $\Ca$ has---again up to a constant---the \emph{minimal} packing $h$-premeasure within $\mathfrak{C}_a$. As with the Hausdorff case, the sharp constant and the exact minimising set are unknown in general; for the middle-third Cantor set, $\mathcal P^s(\Ca) = 4^s \ne \limsup_n n(x_n/n)^s = 1$ where $s=\log2/\log3$, so the elementary bound is not sharp even in this classical case.
\end{proposition}

\begin{theorem}[\cite{HareMendivilZuberman}]
\label{thm:hmz-packing-dichotomy}
Let $h\in\mathcal D$ be subadditive (e.g.\ any $h$ equivalent to a concave function, which includes every associated dimension function $h_a$). The following are equivalent:
\begin{enumerate}[label=(\roman*)]
\item some $E\in\mathfrak{C}_a$ has $\mathcal P^h_0(E)>0$;
\item $\mathcal P^h_0(D_a)=\infty$;
\item $\sum_n h(a_n)=\infty$;
\item some \emph{perfect} $E\in\mathfrak{C}_a$ has $\mathcal P^h(E)=\infty$.
\end{enumerate}
Consequently $\{\mathcal P^h(E) : E\in\mathfrak{C}_a\}$ is either $\{0\}$ or $[0,\infty]$: there is no intermediate regime for the family of packing $h$-measures attained across $\mathfrak{C}_a$.
\end{theorem}

In particular, whenever $\mathcal P^{h_a}_0(\Ca)>0$---which always holds, since $h_a$ is an associated dimension function  of $\Ca$---there exists a rearrangement $E\in\mathfrak{C}_a$ with \emph{infinite} $h_a$-packing measure, even though $D_a$ itself trivially has zero packing measure (being countable) and $\Ca$ realises the minimal packing premeasure in the family.

\begin{theorem}[\cite{HareMendivilZuberman}]
\label{thm:hmz-packing-attain}
Let $h\preccurlyeq h_a$ and $\gamma>0$. There is a perfect set $E\in\mathfrak{C}_a$ with $\mathcal P^h(E)=\gamma$, and the set $\{E\in\mathfrak{C}_a : \mathcal P^h(E)=\gamma\}$ is dense in $(\mathfrak{C}_a,\rho)$.
\end{theorem}

\begin{remark}
Propositions~\ref{prop:hmz-hausdorff-max} and \ref{prop:hmz-packing-min}, taken together, show that the roles of $\Ca$ and the general complementary sets $E\in\mathfrak{C}_a$ are exactly reversed between the Hausdorff and packing measures: $\Ca$ maximises the former and minimises the latter (up to uniform constants), while it is the more disordered rearrangements $E\in\mathfrak{C}_a$ that can achieve arbitrarily small Hausdorff $h$-measure or arbitrarily large---even infinite---packing $h$-measure.
\end{remark}

\subsection{Box dimension}
The box dimensions of sets in $\mathfrak{C}_a$ were carefully studied in \cite[Chapter 3.2]{falconer1997techniques}. There, Falconer shows that the upper and lower box dimensions of sets in $\mathfrak{C}_a$ depend only on the sequence $a$, and not on the rearrangement. That is, given any complementary set, its lower box dimension equals that of every other complementary set, and the same holds for the upper box dimension. Furthermore, for any $E \in \mathfrak{C}_a$,
\[
\liminf_{n \to \infty} \frac{\log n}{-\log a_n} \leq \underline{\dim}_B E \leq \overline{\dim}_B E \leq \limsup_{n \to \infty} \frac{\log n}{-\log a_n},
\]
and the box dimension exists if and only if these limits coincide.

The upper and lower box dimension of $C_a$ (and therefore of any other complementary set in $\mathfrak{C}_a$) are given by 
\[
\overline{\dim}_B C_a =\limsup_{n \to \infty} \frac{\log n}{-\log (x_n/n)},\qquad \underline{\dim}_B C_a =\liminf_{n \to \infty} \frac{\log n}{-\log (x_n/n)}.
\]
Moreover, it was shown in \cite{tricot1994curves} that 
\[
\limsup_{n \to \infty} \frac{\log n}{-\log (x_n/n)}=\limsup_{n \to \infty} \frac{\log n}{-\log a_n}.
\]
This equivalence does not hold for the lower box dimension. In fact, the inequality
\[
\liminf_{n \to \infty} \frac{\log n}{-\log a_n} \leq \liminf_{n \to \infty} \frac{\log n}{-\log (x_n/n)}
\]
can be strict; see, for instance, \cite[Proposition 5]{CabrelliMendivilMolterShonkwiler}.
\subsection{Intermediate dimension}

In order to compute the $\theta$-intermediate dimension of complementary sets,
Angelini and Molter \cite{AngeliniMolter2025} introduce two auxiliary parameters. For a fixed scale $s_n$,
they consider coverings by sets whose diameters lie within the admissible range
determined by $\theta$. Specifically, the aim is to minimise the quantity
\[
\frac{-m \log 2}{\log s_m},
\]
where the index $m$ ranges over all the values satisfying
$s_n^{1/\theta} \leq s_m \leq s_n$; that is, over all admissible diameters $s_m$
that may be used to cover the set at scale $s_n$, in accordance with the
$\theta$-intermediate dimension constraints. Since $0 < s_m < 1$, we have
\[
\frac{-m \log 2}{\log s_m} = \frac{\log 2}{-\log\big(s_m^{1/m}\big)},
\]
so that minimising the left-hand side over the admissible range is equivalent to
minimising $s_m^{1/m}$. This motivates the following definition.

\begin{definition}\label{def:rho}
Let $a = \{a_n\}_{n \in \mathbb{N}}$ be a positive, non-increasing, summable
sequence, let $\theta \in (0,1]$ and let $r \in \mathbb{N}$. We define
$\gamma_{\theta,a}(r)$ and $\rho_{\theta,a}(r)$ as follows:
\[
\gamma_{\theta,a}(r) := \max \left\{ m \in \mathbb{N} : s_m \geq s_r^{1/\theta} \right\},
\]
and
\[
\rho_{\theta,a}(r) := \argmin_{r \leq m \leq \gamma_{\theta,a}(r)+1} \left\{ s_m^{1/m} \right\}.
\]

\end{definition}
The following results are due to Angelini and Molter \cite{AngeliniMolter2025}.
Throughout, $\rho_{\theta,a}$ is the parameter of Definition~\ref{def:rho}, and
we abbreviate
\[
\Lambda_{\theta,n}(a):=\frac{-\rho_{\theta,a}(n)\log 2}{\log s_{\rho_{\theta,a}(n)}} .
\]

\begin{theorem}[The two extremal sets]\label{thm:extremos}
Let $a=\{a_n\}$ be positive, non-increasing and summable.
\begin{enumerate}[label=\textup{(\roman*)}]
\item If $\lim_{n\to\infty}(s_n/s_{n+1})^{1/n}=1$, then for every
$\theta\in(0,1]$,
\[
\overline{\dim}_\theta\, C_a=\limsup_{n\to\infty}\Lambda_{\theta,n}(a),
\qquad
\underline{\dim}_\theta\, C_a=\liminf_{n\to\infty}\Lambda_{\theta,n}(a)
=\dim_H C_a=\underline{\dim}_B\, C_a .
\]
\item For every $\theta\in[0,1]$,
\[
\overline{\dim}_\theta\, D_a
=\frac{\theta\,\overline{\dim}_B\, C_a}{1-(1-\theta)\overline{\dim}_B\, C_a},
\qquad
\underline{\dim}_\theta\, D_a
=\frac{\theta\,\underline{\dim}_B\, C_a}{1-(1-\theta)\underline{\dim}_B\, C_a}.
\]
\end{enumerate}
\end{theorem}

Part~(i) is a genuine formula rather than a mere bound: for $C_a$ the lower
intermediate dimension is constant in $\theta$ and equal to the Hausdorff
dimension, so the interpolation is degenerate from below, while the upper one
varies. Part~(ii) requires no hypothesis on $a$ beyond summability, and the
value depends on $a$ only through the box dimension, which is the same for every
set in $\mathfrak{C}_a$.

\begin{theorem}[Extremality]\label{thm:sandwich}
Under the hypotheses of Theorem~\ref{thm:extremos}\,(i), for every
$\theta\in(0,1]$ and every $E\in\mathfrak{C}_a$,
\[
\overline{\dim}_\theta\, D_a\le\overline{\dim}_\theta\, E
\le\overline{\dim}_\theta\, C_a,
\qquad
\underline{\dim}_\theta\, D_a\le\underline{\dim}_\theta\, E
\le\underline{\dim}_\theta\, C_a .
\]
\end{theorem}

Thus $C_a$ and $D_a$ bracket the whole family, extending to every
$\theta\in(0,1]$ the classical picture at $\theta=0$, where $C_a$ maximises and
$D_a$ minimises the Hausdorff dimension. The lower bounds hold for every
non-increasing $a$, with no regularity assumption.

\begin{definition}\label{def:level-comparable}
A positive, non-increasing, summable sequence $a=\{a_n\}$ is \emph{level
comparable} if
\[
2<\inf_n\frac{s_n}{s_{n+1}}\le\sup_n\frac{s_n}{s_{n+1}}<\infty ,
\]
equivalently, if there are constants $0<\tau\le\lambda<1/2$ with
$\tau\le s_{n+1}/s_n\le\lambda$ for all $n$.
\end{definition}

The right-hand inequality is automatic whenever $a$ is doubling; the content of
the definition is the left-hand one, namely that the ratios $s_{n+1}/s_n$ stay
bounded away from $1/2$. For a central
Cantor set it says exactly that the ratios of dissection are bounded away from
$0$ and from $1/2$, and it gives $s_k\asymp a_{2^k}$.

\begin{theorem}[Full range]\label{thm:rango-theta}
Let $\theta\in(0,1]$ and let $a$ be level comparable. Then
\begin{gather*}
\{\underline{\dim}_\theta\, E : E\in\mathfrak{C}_a\}
=\big[\underline{\dim}_\theta\, D_a,\ \underline{\dim}_\theta\, C_a\big],\\
\{\overline{\dim}_\theta\, E : E\in\mathfrak{C}_a\}
=\big[\overline{\dim}_\theta\, D_a,\ \overline{\dim}_\theta\, C_a\big].
\end{gather*}
\end{theorem}

Every intermediate value is therefore attained: for fixed $\theta$ and any $s$
in the corresponding interval there is a complementary set realising it exactly.
Theorem~\ref{thm:rango-theta} is the $\theta$-intermediate analogue of the
Besicovitch--Taylor theorem, which is the case $\theta=0$.

Theorem~\ref{thm:rango-theta} describes, for each fixed $\theta$, the set of
values attained across $\mathfrak{C}_a$. A finer question is which \emph{functions}
of $\theta$ arise. As recalled in Section~\ref{sec:dimensions}, Banaji and
Rutar \cite{BanajiRutar2022} characterised the functions that occur as the
intermediate dimensions of a bounded subset of $\mathbb{R}^n$; combining their
result with \cite{AngeliniMolter2025} suggests the following: given a fixed
sequence $a$, which functions $\theta\mapsto h(\theta)$ arise as the
intermediate dimension of a complementary set $E\in\mathfrak{C}_a$?

Some partial results in that direction were obtained in
\cite{AngeliniMolter2025}. If $a$ is level comparable, then for every $s>1$
there is a set $E\in\mathfrak{C}_a$ with
\[
\overline{\dim}_\theta\, E=\max\left\{\frac{\overline{\dim}_\theta\, C_a}{s},\
\frac{\theta\,\overline{\dim}_B\, C_a}{1-(1-\theta)\,\overline{\dim}_B\, C_a}\right\}
=:\psi_s(\theta),
\]
so each $\psi_s$ is realised in $\mathfrak{C}_a$. As $s\to1^{+}$ the functions
$\psi_s$ increase to $\overline{\dim}_\theta C_a$, and as $s\to\infty$ they
decrease to $\overline{\dim}_\theta D_a$, so the family $\{\psi_s\}_{s>1}$ sweeps
a one-parameter curve joining the two extremal profiles of
Theorem~\ref{thm:sandwich}. On the other hand, every $E\in\mathfrak{C}_a$
satisfies
\[
\frac{\theta\,\overline{\dim}_B\, C_a}{1-(1-\theta)\,\overline{\dim}_B\, C_a}
\le\overline{\dim}_\theta\, E\le\overline{\dim}_\theta\, C_a .
\]
A one-parameter family of realised profiles inside a two-sided bound is far from
a characterisation, and the following remains open.

\begin{question}\label{q:profiles}
Let $a$ be a level comparable sequence. Describe the structure of the families
of functions
\[
\mathcal{F}^{U}(a):=\big\{\,\theta\mapsto\overline{\dim}_\theta\, E :
E\in\mathfrak{C}_a\,\big\},
\qquad
\mathcal{F}^{L}(a):=\big\{\,\theta\mapsto\underline{\dim}_\theta\, E :
E\in\mathfrak{C}_a\,\big\},
\]
that is, the collections of all functions arising, respectively, as the upper
and lower $\theta$-intermediate dimensions of complementary sets associated
with $a$.
\end{question}
\subsection{Assouad type dimensions}
The Assouad dimension and the lower Assouad dimension of sets in $\mathfrak{C}_a$ were studied in \cite{GarciaHareMendivil2018}. There, the authors prove that, given a non-increasing summable sequence $a = \{a_n\}$,
\[
\dim_A C_a = \limsup_{n \to \infty} \left( \sup_k \frac{n \log 2}{\log(s_k / s_{k+n})} \right),
\]
and
\[
\dim_L C_a = \liminf_{n \to \infty} \left( \inf_k \frac{n \log 2}{\log(s_k / s_{k+n})} \right).
\]
Moreover, it was shown that $\dim_A D_a = 0$ or $1$, and if $E$ is any complementary set in $\mathfrak{C}_a$, then
\[
\dim_A C_a \leq \dim_A E \leq \dim_A D_a,
\]
and
\[
0 = \dim_L D_a \leq \dim_L E \leq \dim_L C_a.
\]

Furthermore, if the sequence $a$ is doubling, then
\[
\{\dim_A E: E \in \mathfrak{C}_a\} = [\dim_A C_a, 1],
\]
and if, instead of doubling, one assumes that $\inf_n s_n/s_{n+1} > 2$, then
\[
\{\dim_L E: E \in \mathfrak{C}_a\} = [0, \dim_L C_a].
\]

The $\Phi$-dimensions introduced in Section~\ref{sec:dimensions} give a
one-parameter refinement of the Assouad-dimension results recalled above. 

\begin{theorem}[{\cite[Prop.~4.2, Thm.~4.3, Cor.~4.7]{GarciaHareMendivil2019b}}]
\label{thm:phi-complementary}
Let $a$ be a positive, non-increasing, summable sequence and let $\Phi$ be an
admissible function. Then, for every $E\in\mathfrak{C}_a$,
\[
\overline{\dim}_\Phi\, E\le\overline{\dim}_\Phi\, D_a
\qquad\text{and}\qquad
\underline{\dim}_\Phi\, E\ge\underline{\dim}_\Phi\, D_a=0 .
\]
If $a$ is level comparable, then in addition
\[
\overline{\dim}_\Phi\, E\ge\overline{\dim}_\Phi\, C_a
\qquad\text{and}\qquad
\underline{\dim}_\Phi\, E\le\underline{\dim}_\Phi\, C_a ,
\]
so that $C_a$ and $D_a$ are extremal in $\mathfrak{C}_a$ for both dimensions. If
moreover $\Phi(x)\to p$ for some $p\in[0,\infty]$, every intermediate value is
attained:
\begin{gather*}
\{\overline{\dim}_\Phi\, E : E\in\mathfrak{C}_a\}
=\big[\overline{\dim}_\Phi\, C_a,\ \overline{\dim}_\Phi\, D_a\big],\\
\{\underline{\dim}_\Phi\, E : E\in\mathfrak{C}_a\}
=\big[0,\ \underline{\dim}_\Phi\, C_a\big].
\end{gather*}
The same conclusions hold for the quasi-Assouad dimensions, with
$\overline{\dim}_\Phi\, C_a$ and $\underline{\dim}_\Phi\, C_a$ replaced by
$\dim_{qA}\,C_a$ and $\dim_{qL}\,C_a$.
\end{theorem}

For the lower $\Phi$-dimension the pattern is the familiar one: the countable
rearrangement $D_a$ is minimal, with value $0$ because it has isolated points,
and the Cantor set $C_a$ is maximal, exactly as for the Hausdorff and lower
Assouad dimensions. For the upper $\Phi$-dimension, however, the roles are
reversed with respect to the Hausdorff dimension case: here $C_a$ realises the \emph{smallest}
value in $\mathfrak{C}_a$ and $D_a$ the largest. 

\begin{remark}\label{rem:Da-equals-one}
If $a$ is doubling and $\Phi(x)\to0$, then $\overline{\dim}_B D_a>0$ and hence
$\overline{\dim}_\Phi D_a\ge\dim_{qA}D_a=1$, so the upper range in
Theorem~\ref{thm:phi-complementary} is the full interval
$[\overline{\dim}_\Phi C_a,1]$. At the other extreme, if $\Phi(x)\to\infty$ then
$\overline{\dim}_\Phi E=\overline{\dim}_B E$, which is the same for every
$E\in\mathfrak{C}_a$, and the interval degenerates to a point.
\end{remark}

Theorem~\ref{thm:phi-complementary} describes the full range of values across
$\mathfrak{C}_a$; one may instead ask for the \emph{typical} value under a natural
random model. Attach to each index $n\in\mathbb{N}$ an independent tag
$\omega_n$, uniformly distributed on $[0,1]$, and order the gaps by increasing
tag. This induces a random total order on $\mathbb{N}$ in which every ordering
of a finite set of indices is equally likely and disjoint sets of indices are
ordered independently, and hence a random element $E_{\boldsymbol\omega}$ of
$\mathfrak{C}_a$. The following dichotomy identifies a sharp threshold, in the
depth of the scales examined, separating two completely different almost sure
behaviours.

\begin{theorem}[{\cite{GarciaHareMendivil2019}}]\label{thm:almost-sure}
Let $a$ be a level comparable sequence and set
$\Psi(x)=\log|\log x|/|\log x|$. Then, almost surely:
\begin{enumerate}[label=\textup{(\roman*)}]
\item if $\Phi(x)\gg\Psi(x)$ as $x\to0$, then
$\overline{\dim}_\Phi \,E_{\boldsymbol\omega}=\overline{\dim}_\Phi\, C_a$ and
$\underline{\dim}_\Phi\, E_{\boldsymbol\omega}=\underline{\dim}_\Phi\, C_a$;
\item if $\Phi(x)\ll\Psi(x)$ as $x\to0$, then
$\overline{\dim}_\Phi\, E_{\boldsymbol\omega}=1$ and
$\underline{\dim}_\Phi\, E_{\boldsymbol\omega}=0$.
\end{enumerate}
\end{theorem}

Since the Assouad dimensions correspond to $\Phi\equiv0$, which falls in case
(ii), and the quasi-Assouad dimensions arise from functions in case (i),
Theorem~\ref{thm:almost-sure} yields the following.

\begin{corollary}\label{cor:almost-sure}
Let $a$ be a level comparable sequence. Then, almost surely, $E=E_{\boldsymbol\omega}$ satisfies
\[
\dim_A E=1,\qquad \dim_L E=0,
\]
while
\[
\dim_{qA}E=\dim_{qA}C_a,\qquad \dim_{qL}E=\dim_{qL}C_a .
\]
\end{corollary}

That is, the almost sure Assouad dimensions of a random rearrangement are always
extremal, whatever $C_a$ may be, whereas the almost sure quasi-Assouad
dimensions always agree with those of $C_a$. In particular, since a set is
uniformly perfect precisely when its lower Assouad dimension is positive, almost
no rearrangement is uniformly perfect. The critical regime $\Phi\asymp\Psi$ is
not covered by Theorem~\ref{thm:almost-sure}.

% =================================================================

% =================================================================

% =================================================================
\section{Dynamically defined Cantor sets}
\label{sec:dynamical}
% =================================================================

The sections above study the whole family $\mathfrak{C}_a$ of complementary sets through the combinatorial and measure-theoretic data of the defining sequence $a$. From now on we restrict to the distinguished  Cantor set $C_a$. A complementary viewpoint, developed in the second half of García's PhD thesis \cite{GarciaThesis2008}, defines Cantor sets \emph{dynamically}, as attractors of iterated function systems (IFS), and asks the converse question: which sequence-defined Cantor sets $\Ca$ arise this way, with what regularity? This dynamical viewpoint is then applied to a genuinely new question---the arithmetic sum $C_p+C_q$ of two members of the ``non-linear'' family $C_p:=C_{\{1/n^p\}}$---via the classical notion of \emph{thickness}.

\subsection{Iterated function systems and hyperbolic attractors}
\label{sec:ifs}

\begin{definition}[IFS and attractor]
Let $D\subset\R$ be closed and nonempty. A finite family $\{f_0,\dots,f_m\}$ of maps $f_i:D\to D$ is an \emph{iterated function system} (IFS) if each $f_i$ is eventually a contraction: there are $r_i<1$ and $n_0\ge1$ with $|f_i^n(x)-f_i^n(y)|\le r_i|x-y|$ for all $x,y\in D$, $n\ge n_0$.
\end{definition}

\begin{theorem}[Hutchinson \cite{Hutchinson1981}]
\label{thm:hutchinson}
Every IFS $\{f_0,\dots,f_m\}$ has a unique nonempty compact \emph{attractor} $K$ satisfying $K=\bigcup_{i=0}^m f_i(K)$, and $f^n(A)\to K$ in the Hausdorff metric on compacta, for every nonempty compact $A\subset D$ with $f_i(A)\subset A$ for all $i$.
\end{theorem}

When the $f_i$ are similarities with ratios $r_i$, the system is \emph{self-similar} and $K$ is a \emph{self-similar set}; the unique $s\ge0$ solving $\sum_i r_i^s=1$ is the \emph{similarity dimension}, and always $\Ha^s(K)\le|K|^s<\infty$, so $\dimH K\le s$. Under the \emph{open set condition} (an open bounded $V$ with $\bigcup_i f_i(V)\subset V$ disjointly), Hutchinson showed $\Ha^s(K)>0$, so $\dimH K = \dim_B K = s$; Bandt--Graf and Schief later showed the open set condition is in fact \emph{equivalent} to $\Ha^s(K)>0$---but not to the packing analogue, since Peres--Simon--Solomyak \cite{PeresSimonSolomyak2000} exhibited self-similar sets on the line with $\Ha^s(K)=0$ yet $0<\mathcal P^s(K)<\infty$ in their common dimension $s$.

\begin{definition}[Regular system]
An IFS $\{f_0,\dots,f_m\}$ on a closed interval $I$ is \emph{regular} if each $f_i\in C^1(I)$ and there are $0<\alpha<1< 1/\beta$, $c\ge1$, $b\ge1$ with, for every $n$ and every word $\omega\in\{0,\dots,m\}^n$: (contractivity) $\alpha^n<|f_\omega'(x)|<c\beta^n$ for all $x\in I$; (bounded distortion) $b^{-1}\le|f_\omega'(x)|/|f_\omega'(y)|\le b$ for all $x,y\in I$, where $f_\omega:=f_{\omega_1}\circ\cdots\circ f_{\omega_n}$.
\end{definition}

Systems with $\eta$-Hölder continuous derivatives, or derivatives satisfying the Dini condition $\int_0^1 m_{f_i'}(t)\,dt/t<\infty$ (with $m_g$ the modulus of continuity), automatically have bounded distortion. For a regular system, the sub-multiplicative quantity $\sum_{\omega\in\{0,\dots,m\}^n}|f_\omega'(x)|^s$ grows like $\exp(n\Upsilon(s))$ for a function $\Upsilon$ independent of $x$; $\Upsilon$ is continuous, strictly decreasing, with $\Upsilon(-\infty)=+\infty$, $\Upsilon(+\infty)=-\infty$ (so it has a unique zero), and equals $\log\sum_i r_i^s$ in the self-similar case. This is \emph{Bowen's equation} $\Upsilon(s)=0$, originally introduced to compute the Hausdorff dimension of quasicircles \cite{Bowen1979}.

\begin{theorem}[\cite{GarciaThesis2008}]
\label{thm:regular-packing}
If $K$ is the attractor of a regular system and $\Upsilon(s)=0$, then $\mathcal P^s_0(K)<\infty$---\emph{unconditionally}, with no separation hypothesis on the system.
\end{theorem}

\begin{theorem}[{\cite[Theorem 3.23]{GarciaThesis2008}}]
\label{thm:regular-osc}
If a regular system satisfies the open set condition, then $0<\Ha^s(K)$ where $\Upsilon(s)=0$; in particular $s=\dimH K=\dim_B K$.
\end{theorem}

Theorems~\ref{thm:regular-packing} and \ref{thm:regular-osc} together give
$0<\mathcal{H}^{s}(K)\le\mathcal{P}^{s}_0(K)<\infty$ under the open set
condition: the attractor of a regular system is an $s$-set, and the same exponent regularises the covering
and the packing process at once. This is the dynamical counterpart of
Theorem~\ref{thm:hset}, where the gauge $h_a$ played the role now played by $s$---with the difference that there a general dimension function was needed,
whereas here regularity of the system forces a power to suffice.

\begin{definition}[Cookie-cutter set]
A regular two-map system $\{f_0,f_1\}$ on $[0,1]$ is \emph{cookie-cutter} (or \emph{hyperbolic}) if it satisfies the strong separation condition
\[
0=f_0(0)<f_0(1)<f_1(0)<f_1(1)=1
\]
together with a strengthening of bounded distortion (uniform convergence of local distortion ratios to $1$ at small scales); the attractor is a \emph{cookie-cutter} (or \emph{hyperbolic}) \emph{Cantor set}. If the defining derivatives are $\eta$-Hölder continuous, the set is $C^{1+\eta}$-hyperbolic.
\end{definition}

\begin{theorem}[\cite{GarciaThesis2008}]
\label{thm:cookiecutter}
If $C$ is the attractor of a cookie-cutter system and $\Upsilon(s)=0$, then $C$ is a Cantor set (compact, perfect, totally disconnected) with $0<\Ha^s(C)\le\mathcal P^s_0(C)<\infty$, and $\mathcal L(C)=0$.
\end{theorem}

Strong separation also gives $C$ a natural symbolic coding by $\Omega^+=\{0,1\}^{\N}$ that is order-preserving for the lexicographic order on $\Omega^+$ (Bedford \cite{Bedford1991}), which will be used below to encode the scale function.

\subsection{The family $C_p$ as a dynamically defined Cantor set}
\label{sec:cp-dynamical}

Recall from Example~\ref{ex:p-cantor} that $C_p=C_a$ for $a_n=n^{-p}$ with
$p>1$, that $\dim_H C_p=1/p$, and that $C_p$ is a $1/p$-set. The main result of
this subsection realises $C_p$ as the attractor of a cookie-cutter system, and
identifies the optimal smoothness of such a realisation.

\begin{theorem}[{\cite[Theorem 1]{Garcia2011}}]\label{thm:cp-hyperbolic}
The set $C_p$ is the attractor of a cookie-cutter system
$\{f_{p,0},f_{p,1}\}$ whose derivatives are $1/p$-Hölder continuous, so $C_p$ is
$C^{1+1/p}$-hyperbolic. Moreover, $1/p$ is the largest Hölder exponent
attainable: no hyperbolic system respecting the step-$k$ interval structure of
$C_p$ has attractor $C_p$ with derivatives of Hölder exponent greater than
$1/p$. In this sense $C_p$ is strictly $C^{1+1/p}$-hyperbolic.
\end{theorem}

Combined with Theorem~\ref{thm:cookiecutter}, this recovers from the dynamical
side the fact recalled above that $C_p$ is a $1/p$-set: hyperbolicity forces
$0<\mathcal{H}^{s}(C_p)\le\mathcal{P}^{s}_0(C_p)<\infty$ with $s=1/p$.

This raises the converse question for \emph{central} Cantor sets $C_a$, i.e.  $a_{2^n+j}=\lambda_n$, $0\le j<2^n$: which decay conditions on $\{\lambda_n\}$ make $C_a$ a dynamically defined set, and with what smoothness? Write $b_n=x_n/n$ as usual and suppose $a$ satisfies the (Besicovitch--Taylor-type) decay condition
\[
\lim_{k\to\infty} 2^k b_{2^k}^s = A, \qquad 0<A<\infty,
\]
which is strictly stronger than merely $0<\varliminf n b_n^s\le\varlimsup n b_n^s<\infty$ (i.e.\ $0<\Ha^s(\Ca)\le\mathcal P^s_0(\Ca)<\infty$, Theorem~\ref{thm:gms-packing}) and is genuinely needed: there are non-increasing sequences satisfying the latter but for which no continuously differentiable system has attractor $C_a$.

\begin{theorem}[\cite{GarciaThesis2008}]
\label{thm:central-c1}
If $a$ is central and satisfies the decay condition above, then $C_a$ is the attractor of a continuously differentiable system $\{f_0,f_1\}$ (with $f_0,f_1$ translates of each other) whose derivative equals the constant $2^{-1/s}$ at every point of $C_a$ itself.
\end{theorem}

The smoothness of this system away from $C_a$ is governed by how fast the derivative's oscillation on the gaps decays. Writing $u_k:=\lambda_{k+1}/\lambda_k - 2^{-1/s}$ for the discrepancy at step $k$:

\begin{theorem}[\cite{GarciaThesis2008}]
\label{thm:central-smoothness}
With $a$, $s$ as in Theorem~\ref{thm:central-c1}: if $\{u_k\}\in\ell^1$, the system has bounded distortion; conversely, if the system's derivatives satisfy the Dini condition, then $\{u_k\}\in\ell^1$. Moreover, the system is strictly $C^{1+\eta}$-hyperbolic if and only if $u_k \sim \lambda_k^{\eta}$.
\end{theorem}

Theorem~\ref{thm:central-smoothness} is a complete dictionary between the decay rate of the sequence and the Hölder smoothness of the dynamically-realising system, and it recovers Theorem~\ref{thm:cp-hyperbolic} for $C_p$ as a special case (with $\eta=1/p$). A characterisation of central sets admitting Hölder-continuously-differentiable realising systems, by different methods, was independently obtained by Bamón--Moreira--Plaza--Vera \cite{BamonMoreiraPlazaVera1997}.

\begin{definition}[Sullivan's scale function]
Let $C$ be a $C^{1+\eta}$-hyperbolic Cantor set with symbolic coding $\Omega^+$; let $\Omega^-$ denote the analogous coding of the ``dual'' Cantor set obtained by composing along the opposite (past) direction. For $\alpha\in\Omega^-$, let $R_n(\alpha)\in\R^3$ record the ratios (left interval, gap, right interval)/(total) at step $n$ around $\alpha$; these converge, as $n\to\infty$, to a limit $R(\alpha)$ that is $\eta$-Hölder continuous. The map $R:\Omega^-\to\R^3$ is the \emph{scale function} of $C$.
\end{definition}

\begin{theorem}[Sullivan \cite{Sullivan1988}, proof in \cite{BedfordFisher1997}]
\label{thm:sullivan}
Two $C^{1+\eta}$-hyperbolic Cantor sets are $C^{1+\eta}$-conjugate if and only if they have the same scale function.
\end{theorem}

That is, the scale function is a complete $C^1$-invariant for hyperbolic Cantor sets---it ``remembers'' the fine geometry of the set at all small scales, and two hyperbolic sets are smoothly equivalent exactly when this invariant coincides.

\begin{theorem}[{\cite[Theorem 2]{Garcia2011}}]
\label{thm:cp-conjugate}
$(C_p,\{f_{p,0},f_{p,1}\})$ is $C^{1+1/p}$-conjugate to the self-similar system
\[
\big(\Lambda_{2^{-p}},\ \{2^{-p}x,\ 2^{-p}x+(1-2^{-p})\}\big),
\]
via an order-preserving diffeomorphism $h_p:[0,1]\to I_p$ with $1/p$-Hölder derivative.
\end{theorem}

Theorem~\ref{thm:cp-conjugate} follows from Theorem~\ref{thm:sullivan} by an explicit computation showing the two sets share the same (constant) scale function. It transports the self-similar (natural) measure $\mu_{2^{-p}}$ on $\Lambda_{2^{-p}}$---the classical \emph{Bernoulli convolution} associated with parameter $2^{-p}$---to the invariant measure $\upsilon_p:=\mu_{2^{-p}}\circ h_p^{-1}$ on $C_p$, which one shows is comparable to $\Ha^{1/p}\llcorner C_p$, since $h_p$ is bi-Lipschitz and Hausdorff measures change by at most bounded factors under bi-Lipschitz maps; we write $\upsilon_p \approx \Ha^{1/p}\llcorner C_p$. This equivalence is the bridge used in the next subsection to transfer results about Bernoulli convolutions to the non-linear family $\{C_p\}$.

\section{Sums of Cantor sets and thickness}
\label{sec:sums}

\subsection{Motivation: Bernoulli convolutions}

The self-similar measure $\mu_r$ on $\Lambda_r$ (with weights $\{1/2,1/2\}$) is the classical \emph{Bernoulli convolution}: the law of $(1-r)\sum_{n\ge0}\varepsilon_n r^n$ for i.i.d.\ fair coin flips $\varepsilon_n\in\{0,1\}$. Every self-similar measure is of \emph{pure type} (singular or absolutely continuous); for $\mu_r$ this dichotomy has been studied intensively since the 1930s.

\begin{theorem}[Erdős, 1939 \cite{Erdos1939}]
\label{thm:erdos}
If $r\in (1/2,1)$ and $1/r$ is a Pisot number (an algebraic integer $>1$ all of whose Galois conjugates have modulus $<1$), then $\mu_r$ is singular; in fact $\widehat{\mu_r}(\xi)\not\to0$ as $\xi\to\infty$.
\end{theorem}

 Whether the set $S_\perp$ of $r\in(1/2,1)$ with $\mu_r$ singular contains anything beyond reciprocals of Pisot numbers is open ; on the positive side, Solomyak \cite{Solomyak1995} showed $\mu_r$ is absolutely continuous for a.e.\ $r\in(1/2,1)$, and Shmerkin \cite{Shmerkin2014} showed $\dimH(S_\perp )=0$.

\subsection{Thickness and Newhouse's Gap Lemma}

A general bound on the dimension of a sum follows from the Lipschitz map $(x,y)\mapsto x+y$: if $\dimH F=\dim_B F$ then $\dimH(E+F)\le\dimH E+\dimH F$ (via the product formula for box-regular factors); since $\dimH C_q=\dim_B C_q$, this gives $\mathcal L(C_p+C_q)=0$ whenever $1/p+1/q<1$.

For the complementary question---when does a sum \emph{contain} an interval---the classical tool is Newhouse's notion of thickness, central to the study of homoclinic bifurcations \cite{PalisTakens1993}.

\begin{definition}[Thickness]
Let $C$ be a Cantor set on the line with a gap $L$. A \emph{bridge} of $L$ is a maximal interval sharing an endpoint with $L$ that meets no gap at least as long as $L$. The \emph{thickness} of $C$ is
\[
\tau(C) := \inf\Big\{ \frac{|B|}{|L|} : (B,L) \text{ is a bridge/gap pair of } C \Big\}.
\]
Thickness is invariant under translation and scaling: $\tau(C+t)=\tau(sC)=\tau(C)$.
\end{definition}

\begin{lemma}[Newhouse's Gap Lemma \cite{newhouse1979abundance}]
\label{lem:newhouse}
If $C_1,C_2$ are Cantor sets on the line with $\tau(C_1)\cdot\tau(C_2)\ge1$, then either $C_1$ is contained in a gap of $C_2$, or $C_2$ is contained in a gap of $C_1$, or $C_1\cap C_2\ne\emptyset$.
\end{lemma}

\begin{proposition}[\cite{PalisTakens1993}]
\label{prop:sum-interval}
If $\tau(C_1)\cdot\tau(C_2)\ge1$, then $C_1+C_2$ is a finite union of closed intervals (so contains an open interval); if in addition no translate of either set lies in a bounded gap of the other, $C_1+C_2$ is a single interval.
\end{proposition}

For example, $\tau(\Lambda_r)=r/(1-2r)$, so $\Lambda_r+\Lambda_r=[0,2]$ once $r\ge1/3$. For $C_p$, a direct computation with the step-interval lengths gives the two-sided bound
\[
\frac{1}{2^p}\cdot\frac{1}{2^p-2} \ \le\ \tau(C_p) \ \le\ \Big(\frac23\Big)^p \frac{1}{2^p-2},
\]
which, via Proposition~\ref{prop:sum-interval}, only yields that $C_p+C_q$ contains an interval under the condition $\frac{1}{2^p(2^p-2)}\cdot\frac1{2^q(2^q-2)}\ge1$---much more restrictive than the condition $\frac1{2^p-2}\cdot\frac1{2^q-2}\ge1$ that one would expect by analogy with the self-similar sets $\Lambda_{2^{-p}}$ and $\Lambda_{2^{-q}}$.

The gap can be closed using \emph{local} thickness $\tau_{\mathrm{loc}}(C,x):=\lim_{k\to\infty}\tau(C\cap I^k_x)$ (the thickness of the sub-Cantor-sets shrinking to $x$), for which the analogous gap-lemma argument only requires the local thickness product to exceed $1$ pointwise:

\begin{theorem}[\cite{Garcia2011}]
\label{thm:local-thickness}
The local thickness of $C_p$ is constant, equal to $\tau_{\mathrm{loc}}(C_p)=1/(2^p-2)$. Consequently, if $\frac{1}{2^p-2}\cdot\frac1{2^q-2}>1$, then $C_p+C_q$ contains an open interval---matching, up to the boundary case, the bound one would guess from the conjugate self-similar sets $\Lambda_{2^{-p}}$.
\end{theorem}

\subsection{Global (almost-everywhere) results}

Beyond the thickness threshold, one can ask for the \emph{exact} Hausdorff dimension or positivity of Lebesgue measure of $C_p+C_q$ across the full range of parameters. Adapting the transversality method of Peres--Solomyak \cite{PeresSolomyak1998} for homogeneous self-similar families---using the $C^{1+1/p}$ conjugacy of Theorem~\ref{thm:cp-conjugate} to control the derivative, in $p$, of the projection $\Gamma_p:\Omega^+\to C_p$---gives:

\begin{theorem}[\cite{Garcia2011}]
\label{thm:global-sum}
For every compact $K\subset\R$: (a) for a.e.\ $p>1$ with $\dimH K+1/p<1$, $\dimH(K+C_p)=\dimH K+1/p$; (b) for a.e.\ $p>1$ with $\dimH K+1/p>1$, $\mathcal L(K+C_p)>0$.
\end{theorem}

\begin{corollary}[\cite{Garcia2011}]
\label{cor:sum-formula}
For fixed $q>1$: for a.e.\ $p>1$, $\dimH(C_p+C_q)=\min\{1/p+1/q,\ 1\}$.
\end{corollary}

Corollary~\ref{cor:sum-formula} is the exact non-linear analogue, for the family $\{C_p\}$, of the Peres--Schmerkin formula $\dimH(\Lambda_r+\Lambda_s)=\min\{\dimH\Lambda_r+\dimH\Lambda_s,1\}$ (valid whenever $\log s/\log r$ is irrational \cite{PeresSchmerkin2007}), obtained here for a full-measure set of the parameter $p$ rather than under an arithmetic condition---reflecting the extra transversality available from the non-linearity of the conjugacy $h_p$.

The absolute continuity half of the Peres--Solomyak method also transfers, via a Frostman-measure argument, to the invariant measures $\upsilon_p$:

\begin{theorem}[\cite{Garcia2011}]
\label{thm:convolution-ac}
Fix $q>1$. For a.e.\ $p>1$ with $1/p+1/q>1$, the convolution $\upsilon_q * \upsilon_p$ is absolutely continuous with density in $L^2(\R)$; equivalently (since $\upsilon_p\approx\Ha^{1/p}\llcorner C_p$), the convolution $\Ha^{1/q}\llcorner C_q * \Ha^{1/p}\llcorner C_p$ is absolutely continuous with $L^2$ density.
\end{theorem}

\begin{remark}
The condition ``for a.e.\ $p$, $q$ fixed'' in Theorem~\ref{thm:convolution-ac} and Corollary~\ref{cor:sum-formula} cannot currently be upgraded to ``for every $p,q$'' or symmetrised in both parameters simultaneously; this is a limitation of the transversality method inherited from \cite{PeresSolomyak1998}, already pointed out in \cite{Garcia2011}. Whether the Fourier transform $\widehat{\upsilon_p}(\xi)$ decays to $0$ as $\xi\to\infty$---the natural non-linear analogue of the Erdős/Salem theory of Section~\ref{sec:sums} above---is entirely open: the non-linearity of the conjugacy $h_p$ destroys the product-of-cosines formula that makes $\widehat{\mu_r}$ tractable in the self-similar case.
\end{remark}

The results of this section connect to a broader circle of ideas on sums of dynamically defined Cantor sets: Palis conjectured that the sum of two $C^{1+\epsilon}$-hyperbolic Cantor sets either has zero Lebesgue measure or contains an interval; this is false in general---Sannami \cite{Sannami1992} constructed a $C^\infty$ Cantor set $A$ such that $A-A$ has positive Lebesgue measure but empty interior, but true generically, in the sense of Moreira--Yoccoz \cite{MoreiraYoccoz2001}, who showed that pairs $(C_1,C_2)$ of hyperbolic sets with $\dimH C_1+\dimH C_2>1$ generically admit a \emph{stable intersection} (a robust transversal intersection of $C_1$ and a translate of $-C_2$), forcing $C_1+C_2$ to contain an interval. Whether $\Lambda_r+\Lambda_s$ has nonempty interior throughout the range $1<\dimH\Lambda_r+\dimH\Lambda_s<1+\delta$ for small $\delta$ remains open.

\begin{openproblem}[\cite{GarciaThesis2008}]
Extend the transversality analysis above to the bi-parametric family $\{C_{p,q}\}_{p>1,q\in\R}$ (sequences $a_n=(\log n)^q/n^p$), using the lexicographic dimension $\dim_{\mathcal F_2}$ induced by $\mathcal F_2:=\{h_{p,q}\}\cup\{x\}\cup\{1\}$ in place of $\dimH$; conjecturally, for fixed $q,q'$ and a.e.\ $(p,p')$ (or fixed $q$ and a.e.\ $p$), $\dim_{\mathcal F_2}(C_{p,q}+C_{p',q'}) = \min\{h_{p,q}\cdot h_{p',q'},\, x\}$ in the sense of the order on $\mathcal F_2$. This would require bi-parametric analogues of Frostman's lemma and of $\beta$-energy adapted to $\dim_{\mathcal F_2}$.
\end{openproblem}

\section{Concluding remarks}
\label{sec:conclusion}

A single pattern runs through the results collected here. For every notion of
size considered, the two canonical members of $\mathfrak{C}_a$ are extremal,
exactly or up to constants, and in most cases where a complete description of
the attained values exists it is a closed interval with endpoints at those two
sets; the packing measures, with their all-or-nothing dichotomy, are the
exception. What changes from one notion to the
next is which of them sits at which end.

For the Hausdorff dimension, and more generally for the Hausdorff $h$-measures,
the maximum is at $C_a$ and the minimum at the countable set $D_a$; the same
holds for the intermediate dimensions at every $\theta\in(0,1]$, and for the
lower Assouad and lower $\Phi$-dimensions. For the upper Assouad and upper
$\Phi$-dimensions the order is reversed, and so it is for the packing
$h$-premeasure: $C_a$ is simultaneously the maximiser of $\mathcal{H}^h$ and the
minimiser of $\mathcal{P}^h_0$ within the family. The reversal is not an
accident of proof technique. Covering-based and global quantities reward the
even distribution of gaps that the Cantor construction produces, while
packing-based and local ones reward the concentration that the decreasing
rearrangement produces, and $C_a$ and $D_a$ are the extreme cases of each.

Against this variability stands the rigidity of the box dimensions, which take
the same value on every member of $\mathfrak{C}_a$. The intermediate dimensions
measure exactly how far a given notion sits from that rigid end: they run from
the full-range behaviour of $\dim_H$ at $\theta=0$ to the constancy of $\dim_B$
at $\theta=1$, and Theorem~\ref{thm:rango-theta} shows that for each
$\theta$ every value between the two extremes is attained, the interval
collapsing to a point only at $\theta=1$.

Several questions remain open. The most substantial, stated as
Question~\ref{q:profiles}, asks not which values but which \emph{functions} of
$\theta$ arise within a fixed family: the ambient answer for arbitrary compact
sets is known from \cite{BanajiRutar2022}, and the profiles $\psi_s$ of
Section~\ref{sec:complementary} give a one-parameter family of realised
examples, but between these lies a genuine gap. For the $\Phi$-dimensions, the
almost sure behaviour in the critical regime $\Phi\asymp\log|\log x|/|\log x|$ is
not covered by the dichotomy of Theorem~\ref{thm:almost-sure}. In the
Hausdorff--packing duality, the sharp constants and the extremal sets are known
on the Hausdorff side only for middle-$\tau$ Cantor sets, and on the packing
side not even for the middle-third Cantor set. Finally, on the dynamical side,
the almost-everywhere results for $C_p+C_q$ cannot at present be made uniform in
both parameters, and the decay of $\widehat{\upsilon_p}$---the non-linear
counterpart of the Erdős--Salem theory for Bernoulli convolutions---is entirely
open.

% =================================================================
\section*{Acknowledgments}
% =================================================================
These notes grew out of three lectures given by the authors at the Instituto de Matemática Aplicada San Luis (Consejo Nacional de Investigaciones Científicas y Técnicas - Universidad Nacional de San Luis), San Luis, Argentina, in
April 2026. We are very grateful to Sergio Favier for the invitation, and  its warm hospitality.
We also thank the participants for their questions and comments, which helped shape this text.


\begin{thebibliography}{99}

\bibitem{AngeliniMolterProjections}
N.~Angelini and U. ~Molter,
\emph{Critical values for intermediate and box dimensions of projections and
other images of compact sets},
J. Fractal Geom. \textbf{12} (2025), no.~3--4, 369--394.

\bibitem{AngeliniMolter2025}
N.~Angelini and U.~Molter,
\emph{Intermediate dimensions of complementary sets},
Proc. Roy. Soc. Edinburgh Sect.~A, published online 2026, 1--19,
\texttt{doi:10.1017/prm.2026.10162}.

\bibitem{AngeliniMolterTejada2026}
N. E. Angelini, U. M. Molter, J. M. Tejada. \emph{Intermediate dimensions of measures: interpolating between Hausdorff and Minkowski dimensions}. J. Math. Anal. Appl. \textbf{555} (2026), no. 1, Paper No. 130039.

\bibitem{BamonMoreiraPlazaVera1997}
R. Bamón, C. G. Moreira, S. Plaza, J. Vera. \emph{Differentiable structures of central Cantor sets}. Ergodic Theory Dynam. Systems \textbf{17} (1997), no. 5, 1027--1042.

\bibitem{Banaji2023}
A. Banaji. \emph{Generalised intermediate dimensions}. Monatsh. Math. \textbf{202} (2023), no. 3, 465--506.

\bibitem{BanajiKolossvary2024}
A.~Banaji and I.~Kolossv\'{a}ry,
\emph{Intermediate dimensions of Bedford--McMullen carpets with applications to
Lipschitz equivalence},
Adv. Math. \textbf{449} (2024), Paper No.~109735, 69~pp.

\bibitem{BanajiRutar2022}
A.~Banaji and A.~Rutar,
\emph{Attainable forms of intermediate dimensions},
Ann. Fenn. Math. \textbf{47} (2022), no.~2, 939--960.

\bibitem{Bedford1991}
T. Bedford. \emph{Applications of dynamical systems theory to fractals---a study of cookie-cutter Cantor sets}. In Fractal Geometry and Analysis (Montreal, 1989), NATO ASI Series C, vol. 346, Kluwer, 1991, 1--44.

\bibitem{BedfordFisher1997}
T. Bedford, A. M. Fisher. \emph{Ratio geometry, rigidity and the scenery process for hyperbolic Cantor sets}. Ergodic Theory Dynam. Systems \textbf{17} (1997), no. 3, 531--564.

\bibitem{BesicovitchTaylor1954}
A. S. Besicovitch, S. J. Taylor. \emph{On the complementary intervals of a linear closed set of zero Lebesgue measure}. J. London Math. Soc. \textbf{29} (1954), 449--459.

\bibitem{Best1939}
E. Best. \emph{A closed dimensionless linear set}. Proc. Edinburgh Math. Soc. (2) \textbf{6} (1939), 105--108.

\bibitem{Bowen1979}
R. Bowen. \emph{Hausdorff dimension of quasicircles}. Inst. Hautes \'Etudes Sci. Publ. Math. \textbf{50} (1979), 11--25.

\bibitem{CabrelliDarjiMolter2011}
C.~Cabrelli, U.~B. Darji and U.~M. Molter,
\emph{Visible and invisible Cantor sets},
in: Excursions in Harmonic Analysis, Vol.~2, 2013.

\bibitem{CabrelliMendivilMolterShonkwiler}
C. Cabrelli, F. Mendivil, U. Molter, R. Shonkwiler. \emph{On the Hausdorff $h$-measure of Cantor sets}. Pacific J. Math. \textbf{217} (2004), no. 1, 45--59.

\bibitem{Davies1971}
R.~O. Davies,
\emph{Sets which are null or non-sigma-finite for every translation invariant measure},
Mathematika \textbf{18} (1971), no.~2, 161--162.

\bibitem{Erdos1939}
P. Erd\H{o}s. \emph{On a family of symmetric Bernoulli convolutions}. Amer. J. Math. \textbf{61} (1939), 974--976.

\bibitem{falconer1997techniques}
K.~J. Falconer.
\newblock \emph{Techniques in Fractal Geometry}. John Wiley \& Sons, Chichester, 1997.

\bibitem{Falconer2014}
K. Falconer. \emph{Fractal Geometry: Mathematical Foundations and Applications}, 3rd ed. Wiley, 2014.

\bibitem{FalconerFraserKempton}
K. Falconer, J. M. Fraser, T. Kempton. \emph{Intermediate dimensions}. Math. Z. \textbf{296} (2020), 813--830.

\bibitem{fraser2020assouad}
J.~M. Fraser,
\emph{Assouad Dimension and Fractal Geometry},
Cambridge Tracts in Mathematics, vol.~222,
Cambridge University Press, Cambridge, 2020.

\bibitem{GarciaThesis2008}
I. A. Garcia. \emph{Análisis fractal de conjuntos de Cantor no lineales}. Tesis Doctoral, Universidad Nacional del Litoral, 2008.

\bibitem{Garcia2011}
I.~Garcia,
\emph{A family of smooth Cantor sets},
Ann. Acad. Sci. Fenn. Math. \textbf{36} (2011), 21--45.

\bibitem{GarciaHareMendivil2018}
I. García, K. Hare, F. Mendivil. \emph{Assouad dimensions of complementary sets}. Proc. Roy. Soc. Edinburgh Sect. A \textbf{150} (2020), 2585--2598.

\bibitem{GarciaHareMendivil2019b}
I.~Garc\'{i}a, K.~E. Hare, and F.~Mendivil,
\emph{Intermediate Assouad-like dimensions},
J. Fractal Geom. \textbf{8} (2021), no.~3, 201--245.

\bibitem{GarciaHareMendivil2019}
I.~Garc\'{i}a, K.~E. Hare, and F.~Mendivil,
\emph{Almost sure Assouad-like dimensions of complementary sets},
Math. Z. \textbf{298} (2021), no.~3--4, 1201--1220.

\bibitem{GarciaMolterScotto2007}
I. Garcia, U. Molter, R. Scotto. \emph{Dimension functions of Cantor sets}. Proc. Amer. Math. Soc. \textbf{135} (2007), no. 10, 3151--3161.

\bibitem{HareMendivilZuberman}
K.~E. Hare, F.~Mendivil and L.~Zuberman,
\emph{The sizes of rearrangements of Cantor sets},
Canad. Math. Bull. \textbf{56} (2013), no.~2, 354--365.

\bibitem{Hutchinson1981}
J. E. Hutchinson. \emph{Fractals and self-similarity}. Indiana Univ. Math. J. \textbf{30} (1981), no. 5, 713--747.

\bibitem{lu2016quasi}
F.~L\"{u} and L.-F. Xi,
\emph{Quasi-Assouad dimension of fractals},
J. Fractal Geom. \textbf{3} (2016), no.~2, 187--215.

\bibitem{MoreiraErazoAngelini}
C.~G. Moreira, H.~Erazo, and N.~Angelini,
\emph{Generalized Hausdorff dimension of irrationals with Lagrange value exactly~$3$},
preprint, 2025, \texttt{arXiv:2511.10328}.

\bibitem{MoreiraYoccoz2001}
C. G. T. de A. Moreira, J.-C. Yoccoz. \emph{Stable intersections of regular Cantor sets with large Hausdorff dimensions}. Ann. of Math. (2) \textbf{154} (2001), no. 1, 45--96.

\bibitem{newhouse1979abundance}
S.~E. Newhouse,
\emph{The abundance of wild hyperbolic sets and non-smooth stable sets for diffeomorphisms},
Publ. Math. Inst. Hautes \'Etudes Sci. \textbf{50} (1979), 101--151.

\bibitem{OlsenRenfro2006}
L.~Olsen and D.~L. Renfro,
\emph{On the exact Hausdorff dimension of the set of Liouville numbers.~II},
Manuscripta Math. \textbf{119} (2006), no.~2, 217--224.

\bibitem{PalisTakens1993}
J. Palis, F. Takens. \emph{Hyperbolicity and Sensitive Chaotic Dynamics at Homoclinic Bifurcations}. Cambridge Studies in Advanced Mathematics, vol. 35, Cambridge University Press, 1993.

\bibitem{PeresSchmerkin2007}
Y.~Peres and P.~Shmerkin,
\emph{Resonance between Cantor sets},
Ergodic Theory Dynam. Systems \textbf{29} (2009), no.~1, 201--221.

\bibitem{PeresSimonSolomyak2000}
Y.~Peres, K.~Simon and B.~Solomyak,
\emph{Self-similar sets of zero Hausdorff measure and positive packing measure},
Israel J. Math. \textbf{117} (2000), 353--379.

\bibitem{PeresSolomyak1998}
Y. Peres, B. Solomyak. \emph{Self-similar measures and intersections of Cantor sets}. Trans. Amer. Math. Soc. \textbf{350} (1998), no. 10, 4065--4087.

\bibitem{Sannami1992}
A.~Sannami,
\emph{An example of a regular Cantor set whose difference set is a Cantor set with positive measure},
Hokkaido Math. J. \textbf{21} (1992), no.~1, 7--24.

\bibitem{Shmerkin2014}
P.~Shmerkin,
\emph{On the exceptional set for absolute continuity of Bernoulli convolutions},
Geom. Funct. Anal. \textbf{24} (2014), no.~3, 946--958.

\bibitem{Solomyak1995}
B. Solomyak. \emph{On the random series $\sum \pm\lambda^n$ (an Erdős problem)}. Ann. of Math. (2) \textbf{142} (1995), no. 3, 611--625.

\bibitem{Sullivan1988}
D. Sullivan. \emph{Differentiable structures on fractal-like sets, determined by intrinsic scaling functions on dual Cantor sets}. In The Mathematical Heritage of Hermann Weyl (Durham, NC, 1987), Proc. Sympos. Pure Math., vol. 48, Amer. Math. Soc., 1988, 15--23.

\bibitem{tricot1994curves}
Claude Tricot.
\newblock {\em Curves and fractal dimension}.
\newblock Springer Science \& Business Media, 1995.

\end{thebibliography}
\end{document}